\documentclass[11pt]{article}

\usepackage[utf8]{inputenc}
\usepackage[T1]{fontenc}

\newcommand{\edge}[1]{\ar@{-}[#1]}

\def\thebibliography#1{\section*{Literature Cited\markboth
 {Literature Cited}{Literature Cited}}\list
 {[\arabic{enumi}]}{\settowidth\labelwidth{[#1]}\leftmargin\labelwidth
 \advance\leftmargin\labelsep
 \usecounter{enumi}}
 \def\newblock{\hskip .11em plus .33em minus -.07em}
 \sloppy
 \sfcode`\.=1000\relax}

\usepackage[all,2cell,matrix,curve,arrow,tips,frame]{xy}
\usepackage{epsfig}
\usepackage{amsmath}
\usepackage{amsfonts}
\usepackage{amsthm}
\usepackage{epic}
\usepackage{eepic}

\newtheorem{Th}{Theorem}[section]
 
\newtheorem{Cy}[Th]{Corollary} 
\newtheorem{Lmm}[Th]{Lemma} 
 
\newtheorem{Defn}[Th]{Definition}

\begin{document}
\begin{center}

 {\LARGE 
 Time Homogeneous Diffusions with a Given Marginal at a Deterministic Time} 
\vspace{.1in}

John M. Noble\\ 
 Institute of Applied Mathematics\\ University of Warsaw\\ ul. Banacha 2\\ 02-097 Warszawa, Poland \\
submitted to Stochastic Processes and Applications: 22nd June 2011, accepted subject to revision 23rd December 2012, revised version submitted 7th May 2012
\end{center}

Key Words: Martingale diffusion, Kre\u in strings, marginal distribution.
 
\begin{abstract}
In this article, it is proved that for any probability law $\mu$ over $\mathbb{R}$ with well defined expected value 
and a given deterministic time $t > 0$, there exists a gap diffusion 
with the prescribed law at the prescribed time.

The method starts by constructing a discrete time process $X$ on a  finite state
space, where $X_\tau$ has law $\mu$, for a geometric time $\tau$, independent of
the diffusion. This argument is developed, using a fixed point theorem, to give conditions for existence for
$\tau$ an independent time with negative binomial distribution. Reducing the
time mesh gives a continuous time diffusion with prescribed law for $\tau$ with
Gamma distribution. Keeping $E[\tau] = t$ fixed, the parameters of the Gamma
distribution are altered, giving the prescribed law for deterministic time. An
approximating sequence establishes the result for arbitrary probability measure
over $\mathbb{R}$.
\end{abstract}

\section{Introduction}

\subsection{The Result and Method of Proof} In this article, it is proved that for any probability measure $\mu$ over $\mathbb{R}$ such that $\int_{\mathbb{R}} |x|\mu(dx) < +\infty$, with expectation denoted $e_0(\mu) = \int_{\mathbb{R}} x\mu(dx)$,   there exists a gap diffusion $X$ in the sense of Kotani Watanabe~\cite{KoWa} and Knight~\cite{Kn},  such that
when $X_0 =  e_0(\mu)$, ${\cal L}(X_t) = \mu$ (where ${\cal L}$ denotes `law'). That is, the probability measure $\mu$ is the law for the random variable $X_t$, the location of the process at a specified
fixed time $t > 0$. The article proves existence; no claims are made about
uniqueness, although the method of proof indicates how a solution may be approximated. 

The proof proceeds in five steps:
\begin{enumerate} 
 \item Discrete time and finite state space are considered; conditions under which a suitable Markov chain with a given distribution when stopped at an independent geometric time are established (theorem~\ref{thgt}, proved in section~\ref{secthgt}). The solution, when it exists is unique and the construction is explicit. 
 
\item This is then extended to establish conditions under which there exists a Markov chain with a given distribution when stopped at an independent negative binomial time. This uses the fact that a negative binomial variable is the sum of independent identically distributed geometric variables and uses a fixed point theorem (theorem~\ref{thnb}, proved in section~\ref{secthnb}). The extension from geometric to negative binomial times, proved in section~\ref{secthnb}, is the crux of the article.

\item The independent time $\tau \sim NB(r,a)$ has expectation $\mathbb{E}[\tau] = \frac{ra}{(1-a)}$. Let $(X^{(a)}_n)_{n \geq 0}$ denote the discrete time chain such that $X_\tau^{(a)}$ has the prescribed terminal distribution. A piecewise constant process  $(Y^{(a,\delta)}_t)_{t \geq 0}$ is constructed, where $Y^{(a,\delta)}_{n\delta} = X_n^{(a)}$. Let $T = \tau \delta$, then $Y^{(a,\delta)}_T$ has the prescribed distribution. The parameter $a \uparrow 1$ and $\delta$ is chosen such that $\frac{ra \delta}{(1-a)} = t$. It follows that $\mathbb{E}[T] = t$ and $\mathbb{V}(T) = t^2\left(\frac{1}{r} + \frac{\delta}{t}\right) \stackrel{\delta \rightarrow 0}{\longrightarrow} \frac{t^2}{r} \stackrel{r \rightarrow +\infty}{\longrightarrow} 0$. For fixed $r$, as $\delta \rightarrow 0$, the distribution of $T$ converges to a $\Gamma(r, \frac{t}{r})$ distribution, while $Y^{(a,\delta)}$ converges to a continuous time Markov process $\tilde{Y}^{(r)}$ on the discrete state space, with the prescribed terminal distribution at the independent gamma distributed time $T$  (theorem~\ref{thgdl1}, proved in subsection~\ref{appthgdlsub1}). 

\item The limit $r \rightarrow +\infty$ is now taken; $T$ converges to $t$ in probability, while $\tilde{Y}^{(r)}$ converges to a continuous time Markov process on a discrete state space with the prescribed terminal distribution at time $t$.
(theorem~\ref{thgdl2}, proved in subsection~\ref{appthgdlsub2}).

\item Finally, arbitrary state space is considered. The target measure $\mu$ is approximated by a sequence of atomised measures $\mu^{(n)}$, where the size of each atom is either approximately the size of the atom of $\mu$ at that point, or  approximately $\frac{1}{2^n}$ if $\mu$ has no atom at that site. The sites of the atoms define a finite state space and the previous arguments show existence of a gap diffusion $Y^{(n)}$ with marginal distribution $\mu^{(n)}$ at a prescribed time $t > 0$. An example shows how continuous time Markov chains on finite state space can be described in the language of gap diffusions. Large $n$ is considered and, by taking subsequences, it is shown that there exists a string measure $m^*$ and a corresponding gap diffusion process, which has marginal distribution $\mu = \lim_{n \rightarrow +\infty} \mu^{(n)}$ at time $t > 0$. This is the content of theorem~\ref{thgdcont}, which is proved in sections~\ref{seckrstr} and~\ref{appthgdcont}, the first of these to set up the background from Kre\u in strings, the second of these to prove the result. 
\end{enumerate}

\subsection{Background} The problem of constructing a gap diffusion with a given law with compact support at an independent exponential time has been discussed fully by Cox, Hobson and Ob\l{}\'{o}j in~\cite{CHO}.
The problem of constructing a martingale diffusion that has law $\mu$ at a fixed
time $t$ has been solved by Jiang and Tao in~\cite{JT} under certain smoothness assumptions.

Recently, Martin Forde in~\cite{F} extended the work of Cox, Hobson and Ob\l{}\'{o}j~\cite{CHO} to provide a process with prescribed joint law for the process at an independent exponential time $\tau$ and its supremum over the time interval $[0,\tau]$. The approach taken by Forde uses the correspondence between the resolvent and the distribution of the process stopped at an exponential time. The argument  requires the probability measure to have a strictly positive density $f$ on the region $R = \{ (y,b) \in {\bf R}^2 | y \leq b, x_0 \leq b < +\infty\}$ where $y$ denotes the process value at the exponential time and $b$ the value of the supremum over the time interval and $x_0 = \int_R y f(y,b) dy db$. The argument is involved, but it should be possible to a) consider the analagous discrete time processes on finite state space stopped at geometric times, b) reduce the time discretisation to obtain continuous time process on finite state space stopped at exponential times  (and thus relax the condition on the density being strictly positive), c)  apply a similar fixed point theorem described here to extend the result to negative binomial times for discrete time steps and gamma times for a continuous time and d) hence take a limit to obtain an existence result for deterministic times on arbitrary state space. The approximating procedure to obtain the general result should be possible. This is a large agenda and the argument by Forde is already involved, but the key ingredients and the role of the resolvent in the proof make the agenda outlined above look possible.   

Since writing this article, the author was made aware (on 19th May 2011) of a
work dated 9th May 2011 submitted, but  unpublished at the time of writing, by Ekström, Hobson,
Janson and Tysk~\cite{EHJT}, that solves the problem in a different way. The article~\cite{EHJT} also approximates the target distribution by atomic measures, but appeals to general results from algebraic topology to conclude existence of a limit. 

In~\cite{Mon}, Monroe constructs a general symmetric stable process with a
prescribed marginal at a fixed time, but does not require that the resulting
process has the martingale property.

It is hoped that the method presented here, although it only proves existence,
may provide the basis of a construction in cases of interest. The key to
existence is showing that there exists a point that satisfies a system of
polynomial equations and there are reasonably efficient numerical methods
available for locating solutions to such systems of equations when they are
known to exist. The discussion in the conclusion indicates the further work
necessary if this is to be developed into a computationally efficient method. 

\subsection{Motivation}

The subject of strong Markov processes generated by Krein-Feller generalised second order differential operators is of great interest in its own right. More specifically, the inverse problem, of computing a string $m$ to give a solution $f$ to the parabolic equation $\frac{\partial f}{\partial t} = \frac{\partial^2 f}{\partial m \partial x}$ with prescribed initial condition at $t = 0$ and prescribed behaviour at $t = T> 0$  is a long standing problem, of interest in mathematical physics.   

In recent years, the interest in the problem has been strongly renewed by  applications to the field of modelling financial markets. Since this is the current driving force for this problem in the applied literature, the financial motivation will be discussed here. 

The general motivating problem within finance is that of automating the pricing and risk management of derivative securities. More specifically, it is the problem of pricing a wide of European style options    given the current market price of the underlying asset and market option quotes of European call options at a range of strikes $K$ for a term $t$ or, more generally, several terms. Here the problem of providing a process that meets a single smile is discussed, but the method could be extended relatively easily to provide a piecewise time homogeneous process that meets given smiles at terms $t_1,\ldots, t_m$. The problem of constructing a process to facilitate option pricing is discussed by Peter Carr in~\cite{C}, who develops a suitable model, known as the {\em local variance gamma model}. This is a problem of practical importance; from listed option prices, the problem of inferring option prices  at non-listed strikes and terms arises both on exchanges and with over-the-counter transactions. The problem  of determining the appropriate inputs for a model so that the output is consistent with a specified set of market prices is known as {\em calibration}.

One of the simplest examples of a calibration procedure is the computation of the implied volatility from the Black Scholes formula. From a single option price, the volatility input for the Black Scholes model is computed, so that pricing is consistent with the given market price. When several prices are given, each with a different maturity, the instantaneous volatility can be considered as a piecewise constant function of time, which jumps at each option maturity. When different strikes each with the same maturity are considered, the implied volatility smile makes it difficult to extend the Black Scholes model in a straightforward way to deal with the set of information. 

Several ways have been suggested to deal with the fact that the implied volatility at a single term is not constant as a function of the strike price and there are many ways to construct a model that is consistent with a given set of arbitrage free market option prices. One approach is found in Rubinstein~\cite{R}, which presents a discrete time model, where the price process is a Markov process on a binomial lattice. A continuous time and state version of Rubinstein's model can be found in Carr and Madan~\cite{CM}. Madan and Yor~\cite{MY} give an alternative way to construct a martingale diffusion that is consistent with a volatility smile. 

The approach of Peter Carr in~\cite{C} is essentially different; the resulting risk neutral process for the price of an asset underlying a set of European options designed to meet a single smile is a time homogeneous process, which is not a diffusion. It is based on a driftless time homogeneous diffusion, which is run on an independent gamma clock. That is, if $X$ denotes the driftless time homogeneous diffusion, then the stock price process $S$ is given by $S_r = X_{\Gamma_r}$ where $\Gamma$ is an independent gamma subordinator. A {\em subordinator} is a one dimensional Lévy process which is increasing almost surely; for a {\em gamma subordinator}, the Lévy process is a gamma process. The gamma clock is normalised so that $\Gamma_t$ has an exponential distribution, where $t$ is the maturity of the options whose prices are given or observed.

This article considers the situation where, for a single fixed term $t$, the European call option, or put option, prices are listed over the whole range of strikes $K$ and shows existence of a risk neutral measure under which the stock price process evolves according to a martingale diffusion where the distribution at time $t$ is that defined by the data. 

\section{Definitions and Results}

\subsection{Markov martingale random walks and gap diffusions}

 The following processes will be used in the article.

\begin{Defn}[Discrete time Markov martingale random walk on a finite state space]
\label{defdtgd} A {\em discrete time Markov martingale random walk} (henceforth referred to as a DMRW) on a finite state space $S =
\{i_1, \ldots, i_M\} \subset \mathbb{R}$ with $i_1 < i_2 < \ldots < i_M$ is a
martingale that is a time homogeneous Markov process with a one step transition
function $P$ with entries $P_{j,k} = \mathbb{P}(X_{t+1} = i_k|X_t = i_j)$ that
satisfy the following conditions: there is a $\underline{q} = (q_1, \ldots, q_M) \in
\{0\} \times [0,1]^{M-2} \times \{0\}$ (taken as a row vector), that is $q_1 = q_M = 0$, such that 

\begin{enumerate}
\item 
\[ P_{j,j} = 1 - q_j, \qquad 0 \leq q_j \leq 1 \qquad j \in \{1,
\ldots, M\}\]
 \item For each $j \in \{2, \ldots, M-1\}$,  
 \[P_{j, {j-1}} = q_j\frac{i_{j+1} - i_j}{i_{j+1} - i_{j-1}}, \qquad
P_{j, {j+1}} = q_j\frac{i_j - i_{j-1}}{i_{j+1} - i_{j-1}}.\]
 \item If $k \not \in \{j-1, j, j+1\}$, then $P_{j,k} = 0$.
 \item For each $j \in \{2,\ldots, M\}$ and each $y \in  \left (i_{j-1},i_j
\right )$, 
 \[\mathbb{P} (X_{0+} = i_j | X_0 = y) = \frac{y - i_{j-1}}{i_j - i_{j-1}} \qquad \mathbb{P}(X_{0+}
= i_{j-1}|X_0 = y) = \frac{i_j - y}{i_j - i_{j-1}}.\]
 \noindent That is, if the initial value $y$ of the process is not in $S$, then
the process immediately jumps (at time $0$) to the  value $\max\{x \in S|x <
y\}$ or the  value $\min\{x \in S | x > y\}$, with probabilities determined to
ensure that the process is a martingale. 
 \end{enumerate}
 \end{Defn}
\paragraph{Notation} The following notation will be used:
\begin{equation}\label{eqalpha} \alpha_{j,j+1} = \frac{i_j -
i_{j-1}}{i_{j+1} - i_{j-1}} \qquad \alpha_{j,j-1} = \frac{i_{j+1} -
i_j}{i_{j+1} - i_{j-1}} \qquad \alpha_{j,k} = 0 \quad k \neq j \pm 1
\end{equation}
\noindent Note that for each $j \in \{1, \ldots, M\}$, $\sum_{k=1}^M P_{jk} = 1$.

\noindent This is a random walk with state space $S$ whose transitions are only  to nearest neighbours in $S$. 

\begin{Defn}[Continuous time Markov martingale random walk on a finite state
space]\label{defctgd} 
A {\em continuous time martingale random walk}(henceforth referred to as a CMRW) on a finite state space $S = \{i_1,
\ldots, i_M\} \subset \mathbb{R}$ with $i_1 < \ldots < i_M$ is a martingale that is
a time homogeneous Markov chain that satisfies the following: there exists a
$\underline{\lambda} = (\lambda_1, \ldots, \lambda_M) \in \{0\} \times \mathbb{R}_+^{M-2} \times \{0\}$ (taken as a row vector with
$\lambda_1 = \lambda_M = 0)$ such that 
\begin{enumerate}
\item \[ \mathbb{P}(X_{s + r} = i_j \quad \forall \quad 0 \leq r \leq t \quad | \quad X_s = i_j) =  e^{-\lambda_j t} \qquad \forall t \geq 0 \qquad j = 1, \ldots, M\]
\item \[ \left\{\begin{array}{l} \lim_{h \rightarrow 0}\frac{1}{h} \mathbb{P}(X_{t + h} = i_{j+1}|X_t = i_j) = \frac{i_j - i_{j-1}}{i_{j+1} - i_{j-1}}\lambda_j \\
\lim_{h \rightarrow 0}\frac{1}{h} \mathbb{P}(X_{t + h} = i_{j-1}|X_t = i_j) =
\frac{i_{j+1} - i_j}{i_{j+1} - i_{j-1}}\lambda_j \qquad   \forall j = 2,
\ldots, M-1\end{array}\right. \]
\item For each $j \in \{2, \ldots, M\}$ and each $y \in (i_{j-1}, i_j)$, 
\[ \mathbb{P}(X_{0+} = i_j | X_0 = y) = \frac{y - i_{j-1}}{i_j - i_{j-1}} \qquad \mathbb{P}(X_{0+}
= i_{j-1}|X_0 = y) = \frac{i_j - y}{i_j - i_{j-1}}.\]
\end{enumerate}
\end{Defn}

\paragraph{Notations} The following notation will be used throughout:
\begin{enumerate} \item 
${\cal B}$ denotes the $\sigma$ algebra of Borel sets.
\item For a measure $\nu$, defined over ${\cal B}(\mathbb{R})$, or the Borel sets of a subset of $\mathbb{R}$, an atom at a point $x \in \mathbb{R}$ is
denoted by $\nu(\{x\})$, where $\{x\}$ denotes the set containing the single point $x \in \mathbb{R}$. 
\end{enumerate}

\begin{Defn}[Continuous time gap diffusion on a continuous state space]\label{defmpsting}  
Let $W(.,x)$ denote a standard Wiener process, with initial condition $W(0,x) = x$. Let $\{\phi^{(x)}(t,z) : t \geq 0, \quad z \in \mathbb{R}\}$ denote the local time of $W(.,x)$ Recall, for example, Revuz and Yor~\cite{RY} chapter 6) that $\phi^{(x)}(.,.)$ is jointly continuous and that $\int_0^s {\bf 1}_A (W(r,x))dr = 2 \int_A \phi^{(x)}(s,z) dz$ for every $A \subseteq {\cal B}(\mathbb{R})$. 

Let $m^*$ denote a measure over $\mathbb{R}$, such that for two points $L_0 < L_1$, possibly $L_0 = -\infty$ and possibly $L_1 = +\infty$, $m^*((a,b)) < +\infty$ for all $L_0 \leq a \leq b \leq L_1$, $-\infty < a < b < +\infty$, $m^*(\{L_0\}) = m^*(\{L_1\}) = +\infty$, $m^*((-\infty,L_0)) = m^*((L_1,+\infty)) = 0$.

Let $T(x,s) = \int_{\mathbb{R}}\phi^{(x)}(s,z)m^*(dz)$ and let $T^{-1}(x,s)$ denote the inverse function of $s \mapsto T(x,s)$. Let

\begin{equation}\label{eqexrep} X(s,x) = W(T^{-1}(x,s),x).
\end{equation}

\noindent Then $X$ defines a strong Markov process on $[L_0, L_1]$ such that $\mathbb{P}(X_{t+s} = L_0 | X_t = L_0)$ for all $t \geq 0$ and all $s \geq 0$ and $\mathbb{P}(X_{t+s} = L_1 | X_t = L_1) = 1$ for all $t \geq 0$ and all $s \geq 0$. 

The process thus defined is the gap diffusion associated with $m^*$. 
 \end{Defn}

\noindent   The definition  is found in Kotani and Watanabe~\cite{KoWa} page 245. 
The terminology {\em gap diffusion} to describe such a process was first introduced by F. Knight. At approximately the same time, S. Kotani and S. Watanabe introduced the terminology {\em generalised diffusion} to discuss the same type of process. The reader is referred to F. Knight~\cite{Kn} and S.Kotani and S. Watanabe~\cite{KoWa} for details.
 
Kotani and Watanabe in~\cite{KoWa} develop the characteristic function, various properties and use it to study properties of generalised diffusions. Knight in~\cite{Kn} considers the local time of the gap diffusion and characterises the set of L\'evy processes that can be obtained by varying the speed measure $m^*$.

\subsection{Results}  
This subsection describes the main results for proving existence of processes with prescribed terminal behaviour. Conditions under which  a DMRW (definition~\ref{defdtgd}) may be constructed, with prescribed distribution when the process is stopped at an independent geometric distribution are given. The construction is shown and the solution, when it exists, is unique. This is the content of theorem~\ref{thgt}.   The quantity ${\cal F}$ defined in definition~\ref{defeff} appears in the explicit formula for the parameters $\underline{q}$ for the DMRW with required terminal distribution at an independent geometric time given by equation (\ref{eqqq}). When $\underline{q} \in \{0\} \times [0,1]^{M-2} \times \{0\}$, equation (\ref{eqqq}) provides the unique solution to the problem; if not, then there is no solution.

Next, the result is extended to show conditions guaranteeing existence of a DMRW with prescribed behaviour when stopped at an independent negative binomial time. This is the content of theorem~\ref{thnb}. This is used to show that there exists a CMRW (definition~\ref{defctgd}) with prescribed behaviour when stopped at an independent time with Gamma distribution. This is the subject of theorem~\ref{thgdl1}. By taking an appropriate limit, so that the coefficient of variation goes to zero, a CMRW with prescribed behaviour at a deterministic time is obtained. This is the subject of theorem~\ref{thgdl2}. By considering the cumulative distribution function in the general case as the limit of cumulative distribution functions of variables with finite state space, the general result, existence of a gap diffusion with a prescribed law at a fixed time is obtained. 

\begin{Defn}\label{defeff} Let  $\underline{p} = (p_1, \ldots, p_M)$ (taken as a row vector) be a probability
mass function, that is  $p_j \geq 0$ for each $j \in \{1, \ldots, M\}$ and $\sum_{j=1}^M p_j = 1$. Let $S = \{i_1, \ldots, i_M\} \subset \mathbb{R}$, $i_1 <
\ldots < i_M$ be the support of the probability mass function and $e_0(\underline{p})$, or simply $e_0$ when it is clear which expectation is meant, denote its expected value;
\[ e_0(\underline{p}) = e_0 = \sum_{j=1}^M i_j p_j.\]

\noindent Let $l$ denote the coefficient such that $i_{l-1} < e_0 \leq i_l$. Set  

\begin{equation}\label{eqeldef} 
{\cal L}(\underline{p},j) = \left\{\begin{array}{ll}\frac{(i_{j+1} - i_{j-1})}{(i_{j+1} -
i_j)(i_j - i_{j-1})} \left( \sum_{k=1}^{j-1} (i_j - i_k)p_k\right ) & 2 \leq
j \leq l-1  \\
\frac{(i_{j+1} - i_{j-1})}{(i_{j+1} - i_j)(i_j - i_{j-1}) } \left(
\sum_{k=j+1}^{M} (i_k - i_j)p_k \right ) & l \leq j \leq M-1\\
0 & j = 1 \quad \mbox{or} \quad M
\end{array} \right.
\end{equation} 

\noindent and 

\begin{equation}\label{eqefdef} 
{\cal F}(\underline{p},j) = \left\{\begin{array}{ll}\frac{{\cal L}(\underline{p},j)}{p_j}   & j \in \{1,
\ldots, M\} \qquad \mbox{when} \qquad {\cal L}(\underline{p},j) > 0  \\
0 & j \in \{1, \ldots, M\} \qquad \mbox{when} \qquad {\cal L}(\underline{p},j) = 0
\end{array} \right.
\end{equation} 
\end{Defn}

\paragraph{Notation} Throughout, $e_0$ will be used to denote the expected value of a probability distribution. If $\underline{p} = (p_1, \ldots, p_M)$ represents a probability mass function over a finite set $\{i_1, \ldots, i_M\}$, then the notation $e_0$ or $e_0(\underline{p})$ will be used to denote $e_0(\underline{p}) = \sum_{j=1}^M i_j p_j$. If $\mu$ is a probability measure over ${\bf R}$ with a well defined expected value, then the notation $e_0$ or $e_0(\mu)$ will be used to denote $e_0(\mu) = \int x \mu(dx)$. \vspace{5mm}

\noindent The following theorem is a discrete version of the first approach to the problem of finding a process with prescribed terminal distribution at an exponential time by Cox, Hobson and Ob\l{}\'{o}j~\cite{CHO}.  Exponential time is replaced by geometric time, which is its discrete analogue, and the argument is similar to the use of the resolvent in~\cite{CHO}. In the discrete setting, the argument is similar and boils down to showing that there exists a solution to a system of linear equations. Exploiting the idea that a diffusion at an  exponential time could be computed explicitly by considering the resolvent (the approach presented below for geometric times) appeared earlier, in Peter Carr~\cite{C}.
 
 \begin{Th}[DMRW process at geometric time]\label{thgt} Let $\tau$ denote a random time with
probability function
 \[ p_\tau(k) = \mathbb{P}(\tau = k) = \left\{\begin{array}{ll} (1-a)a^k & k =
0,1,2,\ldots \\
0 & \mbox{otherwise} 
                                     \end{array}\right. \]

\noindent That is $\tau \sim Ge(a)$ (geometric distribution with parameter $a$). Let $\underline{p} = (p_1, \ldots, p_M)$ be a probability mass function, satisfying $\min_{j \in \{1,\ldots, M\}}p_j > 0$. Let  $S = \{i_1, \ldots, i_M\} \subset \mathbb{R}$, $i_1 < \ldots < i_M$ and suppose that $\underline{p}$, taken as a probability mass function over $S$ has expected value $e_0 = e_0 (\underline{p} ) = \sum_{j=1}^M i_j p_j$. Let ${\cal F}$ be defined by
equation (\ref{eqefdef}). Then there is a DMRW $X$  
(definition~\ref{defdtgd}) with state space $S$ and one step transition matrix $P$ as in the
definition, where $\tau$ is independent of $X$, such that 
for $l$ satisfying $i_{l-1} < e_0 \leq
i_l$, $i_{l-1}, i_l \in S$,   
 \[\frac{e_0 - i_{l-1}}{i_l - i_{l-1}}\mathbb{P}(X_\tau = i_j | X_0 = i_l) + \frac{i_l -
e_0}{i_l - i_{l-1}}\mathbb{P}(X_\tau = i_j | X_0 = i_{l-1}) = p_j \] 
 
\noindent if and only if

\begin{equation}\label{eqqbd} \max_{k \in \{1, \ldots, M\}} \left(\frac{1}{a} -
1\right) {\cal F}(\underline{p},k) \leq 1. 
\end{equation}

\noindent  The components of vector $\underline{q}$
are given by

\begin{equation}\label{eqqq} \left\{\begin{array}{l} q_j = \left(\frac{1}{a}
- 1\right) {\cal F}(\underline{p},j) \qquad j \in \{2, \ldots, M - 1\} \\ q_1 = q_M = 0.
\end{array}\right. \end{equation}
\end{Th}

\paragraph{Proof of theorem~\ref{thgt}} This is the subject of section~\ref{secthgt}. \qed \vspace{5mm}

\noindent This is then extended, using a fixed point theorem (theorem~\ref{thfpt}), to negative binomial times. A negative binomial time may be regarded as the sum of independent identically distributed geometric times. For a $NB(r,a)$ time (negative binomial, sum of $r$ independent geometric times, each with parameter $a$), this boils down to showing that there is a solution to a system of polynomial equations, each of degree $r$. The fixed point theorem is required to give conditions under which there exists a solution.  
  
\begin{Th} [DMRW process at negative binomial time] \label{thnb} 
Let $S = \{i_1, \ldots, i_M\} \subset \mathbb{R}$. Let $\underline{p} = (p_1, \ldots, p_M)$ be a probability mass
function satisfying $\min_{j \in \{1, \ldots, M\}} p_j > 0$. Suppose that $\underline{p}$ is a probability mass function over  $S$, with expected value  
$e_0 = e_0(\underline{p}) = \sum_{j=1}^M i_j p_j$. Let $\tau$ be a random time with
probability function 

\begin{equation}\label{eqnbtime} p_\tau(k) = \mathbb{P}(\tau = k)
=\left\{\begin{array}{ll} \left(\begin{array}{c} r+k-1 \\r-1 \end{array}\right)
(1-a)^r a^k & k = 0,1,2, \ldots \\ 0 & \mbox{otherwise}
\end{array}\right.\end{equation}

\noindent That is, $\tau \sim NB(r,a)$ (negative binomial with parameters $r$ and $a$). Let $q_1 = q_M = 0$. Then there exists an $a_0 \in [0,1)$ such
that for each $a \in [a_0,1)$ there exists a $\underline{q} = (q_1, \ldots, q_M) \in
\{0\} \times (0,1]^{M-2} \times \{0\}$ such that $P(\underline{q})$ defined in definition~\ref{defdtgd} is the one
step transition matrix for a Markov chain $X$ with state space $S$, independent of $\tau$, such that

\[ \mathbb{P}(X_\tau = i_j |X_0 = e_0) = p_j \qquad j = 1, \ldots, M.\]

\end{Th}

\paragraph{Proof of theorem~\ref{thnb}} The fixed point theorem to enable the result of theorem~\ref{thgt} to negative binomial times, and the proof of theorem~\ref{thnb} are the subject of section~\ref{secthnb}. \qed

\paragraph{Aside} Note that

\[ \mathbb{E}[\tau ] = \frac{ar}{1-a} = \mu \qquad \mathbb{V}(\tau) = \frac{ar}{(1-a)^2}
= \mu^2\left(\frac{1}{r} + \frac{1}{\mu}\right)  \] 

\noindent where $\mu$ is used to denote $\mathbb{E}[\tau]$ and $\mathbb{V}$ denotes variance. The idea is to consider time
steps of length $\delta$, where $\delta \rightarrow 0$, with a terminal time $T
= \delta \tau$, $\mathbb{E}[T] = t$. Then

\[ \mathbb{E}[T] = t, \qquad \mathbb{V}(T) = t^2\left(\frac{1}{r} +
\frac{\delta}{t}\right) \stackrel{\delta \rightarrow 0}{\longrightarrow}
\frac{t^2}{r} \stackrel{r \rightarrow +\infty}{\longrightarrow} 0.\] \qed

\noindent The result of theorem~\ref{thnb} may therefore be extended to show existence of a CMRW (definition~\ref{defctgd}) with a prescribed distribution at a deterministic time $t > 0$. Formalising the argument outlined above is the subject of theorems~\ref{thgdl1} and~\ref{thgdl2}; theorem~\ref{thgdl1} gives existence of a CMRW with prescribed distribution at an independent gamma time, while theorem~\ref{thgdl2} alters the parameters, keeping the expected value fixed, to obtain a CMRW with a prescribed distribution at a deterministic time $t > 0$. 
 
\begin{Th}[CMRW process at gamma time] \label{thgdl1}
Let $S = \{i_1, \ldots, i_M\} \subset \mathbb{R}$ and let $\underline{p} = (p_1, \ldots p_M)$ (taken as a row vector) be a probability mass
function over $S$, where $\min_{j \in \{1,\ldots, M\}}p_j > 0$. That is,
$\sum_{j=1}^M p_j = 1$ and zero probability is assigned to $x \not \in S$. Let $e_0 = e_0(\underline{p}) =
\sum_{j=1}^M i_j p_j$. Let $\tau$ be a random time with probability density function 
\[ f_\tau (x) = \frac{1}{\Gamma(r)} \left(\frac{r}{t} \right)^r x^{r-1} e^{-rx/t} \qquad x \geq 0.\]
That is, $ \tau \sim \Gamma (r, \frac{t}{r})$, gamma distribution with parameters $r$ and $\frac{t}{r}$ where $\Gamma$ denotes the Euler Gamma function, which for integer $r \geq 1$ is $\Gamma(r) = (r-1)!$. For any integer $r \geq 1$ and any $t > 0$, $t \in {\bf R}_+$, there exists a CMRW process $X$ on $S$  satisfying definition~\ref{defctgd} such that 
\[ \mathbb{P}(X_\tau = i_j | X_0 = e_0) = p_j \qquad j = 1, \ldots, M.\]
\end{Th}

\paragraph{Proof of theorem~\ref{thgdl1}} This is the subject of   subsection ~\ref{appthgdlsub1}. \qed 

\begin{Th}[CMRW process at a fixed time] \label{thgdl2}
Let $S = \{i_1, \ldots, i_M\} \subset \mathbb{R}$ and let $\underline{p} = (p_1, \ldots p_M)$ (taken as a row vector) be a probability mass
function satisfying $\min_{j \in \{1,\ldots, M\}}p_j > 0$.  Suppose $\underline{p}$ is a probability mass function over $S$ and let $e_0 = e_0(\underline{p}) =
\sum_{j=1}^M i_j p_j$ denote its expected value. For any specified $t > 0$, there exists a CMRW process $X$ with state space $S$  satisfying definition~\ref{defctgd} such that 
\[ \mathbb{P}(X_t = i_j | X_0 = e_0) = p_j \qquad j = 1, \ldots, M.\]
\end{Th}

\paragraph{Proof of theorem~\ref{thgdl2}} This is the subject of subsection~\ref{appthgdlsub2}. \qed \vspace{5mm}

\noindent The final, and largest step is to take a limit and go from atomised probability measures with finite numbers of atoms to arbitrary probability measures over $\mathbb{R}$. This is the subject of theorem~\ref{thgdcont}. The proof of this theorem requires section~\ref{seckrstr} as preparatory material and  section~\ref{appthgdcont} to prove the convergence.  While the steps are routine, a substantial quantity of analysis is necessary to take the limit. Lemma~\ref{lmmloctimeconv} of section~\ref{appthgdcont} uses the characterisation of Markov processes in terms of a time changed Wiener process, where the time change is given in terms of the local time and the string measure to show that if the strings converge appropriately, then the processes converge in a suitable sense. The task is therefore to prove that, given a suitable approximating sequence of atomised measures approximating the target measure, there is a convergent subsequence of strings with a well defined limit.  Convergence is then considered in two parts; firstly, convergence of part of the strings to a measure that is absolutely continuous and convergence to the atomic part, corresponding to the atoms in the target measure.   

\begin{Th}[Gap Diffusion]\label{thgdcont} For any probability distribution
function $\mu$ defined on $\mathbb{R}$ such that $\int_{-\infty}^\infty |x|\mu(dx) < +\infty$, let $e_0 = \int_{-\infty}^{\infty} x\mu(dx)$ denote its expectation. Let $t$
be a fixed time $t > 0$. Then there exists a string measure $m^*$ over $\mathbb{R}$ such
that $\frac{d^2}{dm^* dx}$ is the infinitesimal generator of a gap diffusion $X$
where 

\[\mathbb{P}(X_t \leq x | X_0 = e_0) = \mu((-\infty, x]) \qquad \forall x \in (-\infty, +\infty).\]
\end{Th}

\paragraph{Proof of theorem~\ref{thgdcont}} This is the subject of section~\ref{seckrstr} to define the machinery and section~\ref{appthgdcont} to prove the convergence. \qed

\section{Proof of theorem~\ref{thgt}}\label{secthgt} This section is devoted to the proof of theorem~\ref{thgt}, giving conditions under which there exists a DMRW (definition~\ref{defdtgd}) with a given distribution at an independent geometric time, together with an explicit formula. The explicit formula   is of crucial importance in the proof of theorem~\ref{thnb} where the existence result is extended to negative binomial times.  

Firstly, if $\tau \sim Ge(a)$ (Geometric with paramter $a$) independent of $X$, then 
\[ \mathbb{P}(X_\tau = i_k | X_0 = i_j) = \sum_{m = 0}^\infty \mathbb{P}(X_m = i_k | X_0 = i_j) \mathbb{P}(\tau = m) = (1-a) \sum_{m = 0}^\infty (aP)^m_{j,k}  = (1-a) (I - aP(\underline{q}))^{-1}_{jk},\]

\noindent where $(aP)^m = a^m P^m$ and $P^m$ is taken in the sense of multiplication of matrices, so that if the process has initial distribution $\underline{v}$ at time $0+$ and terminal distribution $\underline{p}$ at time $\tau$, then $\underline{p} = (1-a)\underline{v} (I - aP(\underline{q}))^{-1}$. 

The central part of the proof of theorem~\ref{thgt} is lemma~\ref{lmmone}. Lemma~\ref{lmmone}  gives an explicit computation for the parameters $\underline{q}$ such that $\underline{p} (I - aP(\underline{q})) = (1-a) \underline{v}$ for a specified measure $\underline{p}$ and a vector $\underline{v}$ specifying the probability distribution at time $0+$; if $e_0(\underline{p}) \not \in S$, the process jumps into space $S$ at time $0$ according to the rules described in definition~\ref{defdtgd}, which gives the vector $\underline{v}$. The system of equations is linear and has an explicit solution.  The proof of theorem~\ref{thgt} is then a simple corollary.

The quantity ${\cal F}$ of equation (\ref{eqefdef}) definition~\ref{defeff} is of crucial importance in the whole construction; it turns out (lemma~\ref{lmmone}) that it is the explicit formula for the parameters $\underline{q}$ for the DMRW with required terminal distribution at an independent geometric time given by equation (\ref{eqqq}).

\begin{Lmm} \label{lmmone}
Let  $\underline{p} = (p_1, \ldots, p_M)$ denote a row vector such that $\min_j p_j > 0$ and $\sum_{k=1}^M p_j = 1$. Let ${\cal L}$ and ${\cal F}$ be defined as in definition~\ref{defeff}, $\alpha_{j,k}$ by equation (\ref{eqalpha}).     Let $e_0 = e_0(\underline{p}) = \sum_{j=1}^M i_j p_j$. Then  there is a unique row vector $\underline{q} = (q_1, \ldots, q_M) \in  \{0\} \times \mathbb{R}_+^{M-2} \times \{0\}$  and  matrix $P(\underline{q})$
defined by 
\begin{equation}\label{eqPdef} \left\{ \begin{array}{l} P_{j,j}(\underline{q}) = 1 - q_j \qquad j = 1,\ldots, M \\
 P_{j,{j-1}}(\underline{q}) = \alpha_{j,j-1}q_j, \qquad j = 2, \ldots, M \\ 
  P_{ j, {j+1}}(\underline{q}) = \alpha_{j,j+1}q_j \qquad j = 1, \ldots, M-1\\ 
 P_{ j, k} = 0 \qquad k \not \in \{j-1, j, j+1\} \qquad (j,k) \in \{1,
\ldots, M\}^2 \end{array}\right.
\end{equation}

\noindent such that for
each $a \in (0,1)$,   if $i_{l-1} < e_0 \leq i_l$, then $P$ satisfies

\begin{equation}\label{eqmq2} (\underline{p}(I-a P(\underline{q})))_j = \left\{ \begin{array}{ll} 0 & j
\neq l-1, l \\  \frac{e_0 - i_{l-1}}{i_l - i_{l-1}}(1 - a) & j = l \\ \frac{i_l
- e_0}{i_l - i_{l-1}}(1-a) & j = l-1. \end{array}\right. \end{equation}

\noindent The parameters $\underline{q} = (q_1, \ldots, q_M)$ satisfy 

\begin{equation}\label{eqlmmonesol} q_j =  \left(\frac{1}{a} - 1\right)  {\cal F}(\underline{p},j) \qquad j \in \{1, \ldots,
M\}.\end{equation}
\end{Lmm}

\paragraph{Proof of lemma~\ref{lmmone}}  
\noindent From the definition of $P(\underline{q})$ in equation (\ref{eqPdef}), equation (\ref{eqmq2}) may be written as

\begin{equation}\label{eqmq3} \underline{p}{\cal M} = (1-a) \underline{v}
\end{equation}

\noindent where 

\begin{equation}\label{eqemm} \left\{\begin{array}{ll} {\cal M}_{jj} = 1 - a + aq_j & j = 1, \ldots, M \\
 {\cal M}_{j,j+1} =    -aq_j \alpha_{j, {j+1}} & j = 1, \ldots, M-1 \\  {\cal M}_{j, {j-1}} = -aq_j \alpha_{ j,  {j-1}} & j = 2,\ldots, M \\ {\cal M}_{j,k} = 0 & k \not\in \{j-1,j,j+1\}\end{array} \right. \end{equation}
 
 \noindent and
 
 \[ v_j = \left\{ \begin{array}{ll} 0 & j
\neq l-1, l \\  \frac{e_0 - i_{l-1}}{i_l - i_{l-1}}  & j = l \\ \frac{i_l
- e_0}{i_l - i_{l-1}}  & j = l-1. \end{array}\right.\]

\noindent Equation (\ref{eqmq3}) gives a linear system of $M$ equations and $M-2$
unknowns $(q_2, \ldots, q_{M-1})$, but 
\[ \sum_j (\underline{p}(I - aP(\underline{q}))_j = \sum_{jk}p_k(I - aP(\underline{q}))_{kj} = \sum_k p_k (1 - a + aq_k - aq_k \alpha_{k,k-1} - a q_k \alpha_{k,k+1}) = (1-a),\]
\noindent while $ \sum_j v_j = 1$. Furthermore,

\begin{eqnarray*} \lefteqn{\sum_j i_j (\underline{p}(I - aP(\underline{q}))_j}\\&& = \sum_{jk}i_j p_k(I - aP)_{kj} = \sum_k p_k (i_k(1 - a) + a i_k q_k - i_{k-1} aq_k \alpha_{k,k-1} - a i_{k+1} q_k \alpha_{k,k+1})\\&& = (1-a)\sum_k i_k p_k + a \sum_{k=2}^{M-1} p_k q_k \left (1 - i_{k-1}\frac{i_{k+1} - i_k}{i_{k+1} - i_{k-1}} - i_{k+1}\frac{i_{k} - i_{k-1}}{i_{k+1} - i_{k-1}}  \right )\\&& = (1-a)e_0 
\end{eqnarray*}

\noindent while 
\[ \sum_j i_j v_j = i_l \frac{e_0 - i_{l-1}}{i_l - i_{l-1}}  + i_{l-1}\frac{i_l
- e_0}{i_l - i_{l-1}} =  e_0.\] 

\noindent It follows that there are at most $M-2$ linearly independent equations with $M-2$ unknowns and hence at most one solution. It is now  shown that $q_j:j=2, \ldots, M-1$ given by equation (\ref{eqlmmonesol}) provides a solution.\vspace{5mm}

\noindent Consider $q_1 = q_M = 0$ and $q_2, \ldots, q_{M-1}$ defined by equation (\ref{eqlmmonesol}). For $3 \leq j \leq l-2$,  

\begin{eqnarray*}
 \lefteqn{\sum_k p_k (I - aP(\underline{q}))_{kj}}\\&& = (1-a) p_j + (1-a)  {\cal L}(\underline{p},j) - (1-a){\cal L}(\underline{p}, j-1) - (1-a){\cal L}(\underline{p}, j+1) \\
 && = p_j (1-a) + (1-a)\frac{(i_{j+1} - i_{j-1})}{(i_{j+1} - i_j)(i_j - i_{j-1})}\sum_{k=1}^{j-1} (i_j - i_k)p_k\\&& \hspace{5mm}
 - (1-a)\frac{1}{i_j - i_{j-1}} \sum_{k=1}^{j-2}(i_{j-1} - i_k)p_k - (1-a)\frac{1}{i_{j+1} - i_j}\sum_{k=1}^j (i_{j+1} - i_k)p_k 
 \\&& = (1-a) \sum_{k=1}^{j-2} \left \{\frac{(i_{j+1} - i_{j-1})(i_j - i_k)}{(i_{j+1} - i_j)(i_j - i_{j-1})} - \frac{i_{j-1} - i_k}{i_j - i_{j-1}} - \frac{i_{j+1} - i_k}{i_{j+1} - i_j} \right \} p_k\\&& \hspace{5mm}
 + (1-a)\left\{\frac{(i_{j+1} - i_{j-1})(i_j - i_{j-1})}{(i_{j+1} - i_j)(i_j - i_{j-1})} - \frac{i_{j+1} - i_{j-1}}{i_{j+1} - i_j}\right\}p_{j-1}\\&& = 0.
\end{eqnarray*}

\noindent For $2 = j \leq l-2$,

\begin{eqnarray*}
  \sum_k p_k (I - aP)_{k2} &=& (1-a) p_2 + (1-a)  {\cal L}(\underline{p},2)   - (1-a){\cal L}(\underline{p}, 3) \\
 &=& p_2 (1-a) + (1-a)\frac{(i_{3} - i_{1})}{(i_{3} - i_2)(i_2 - i_{1})}  (i_2 - i_1)p_1 
  - (1-a)\frac{1}{i_{3} - i_2}\sum_{k=1}^2 (i_{3} - i_k)p_k 
 \\&=&  (1-a)\left\{\frac{(i_{3} - i_{1})(i_2 - i_{1})}{(i_{3} - i_2)(i_2 - i_{1})} - \frac{i_{3} - i_{1}}{i_{3} - i_2}\right\}p_{1} = 0.
\end{eqnarray*}

\noindent For $l+1 \leq j \leq M-2$, 

\begin{eqnarray*}
 \lefteqn{\sum_k p_k (I - aP)_{kj}}\\&& = (1-a) p_j + (1-a)  {\cal L}(\underline{p},j) - (1-a){\cal L}(\underline{p}, j-1) - (1-a){\cal L}(\underline{p}, j+1) \\
 && = p_j (1-a) + (1-a)\frac{(i_{j+1} - i_{j-1})}{(i_{j+1} - i_j)(i_j - i_{j-1})}\sum_{k=j+1}^{M} (i_k - i_j)p_k\\&& \hspace{5mm}
 - (1-a)\frac{1}{i_j - i_{j-1}} \sum_{k=j}^{M}(i_{k} - i_{j-1})p_k - (1-a)\frac{1}{i_{j+1} - i_j}\sum_{k=j+2}^M (i_{k} - i_{j+1})p_k 
 \\&& = (1-a) \sum_{k=j+2}^{M} \left \{\frac{(i_{j+1} - i_{j-1})(i_k - i_j)}{(i_{j+1} - i_j)(i_j - i_{j-1})} - \frac{i_{k} - i_{j-1}}{i_j - i_{j-1}} - \frac{i_{k} - i_{j+1}}{i_{j+1} - i_j} \right \} p_k\\&& \hspace{5mm}
 + (1-a)\left\{\frac{(i_{j+1} - i_{j-1})(i_{j+1} - i_j)}{(i_{j+1} - i_j)(i_j - i_{j-1})} - \frac{i_{j+1} - i_{j-1}}{i_{j} - i_{j-1}}\right\}p_{j+1}\\&& = 0.
\end{eqnarray*}

\noindent For $M-1 = j \geq l+1$, 

\begin{eqnarray*}
  \sum_k p_k (I - aP)_{k,M-1} &=& (1-a) p_{M-1} + (1-a)  {\cal L}(\underline{p},M-1)   - (1-a){\cal L}(\underline{p}, M-2) \\
 &=& p_{M-1} (1-a) + (1-a)\frac{(i_{M} - i_{M-2})(i_M - i_{M-1})}{(i_{M-1} - i_{M-2})(i_{M} - i_{M-1}) } p_M \\&& \hspace{5mm}  
  - (1-a)\frac{1}{i_{M-1} - i_{M-2}}\sum_{k=M-1}^M (i_{k} - i_{M-2})p_k 
 \\&=&  (1-a)\left\{\frac{(i_{M} - i_{M-2})(i_M - i_{M-1})}{(i_M - i_{M-1})(i_{M-1} - i_{M-2})} - \frac{i_{M} - i_{M-2}}{i_{M-1} - i_{M-2}}\right\}p_{M} = 0.
\end{eqnarray*}

\noindent For $j = l-1$,  $\sum_{k=1}^M i_k p_k = e_0$ gives $\sum_{k= l+1}^M i_k p_k = e_0 - \sum_{k=1}^l i_k p_k$  and $\sum_{k=1}^M p_k = 1$ gives $\sum_{k=l+1}^M p_k = 1 - \sum_{k=1}^l p_k$. Using this,

\begin{eqnarray*}\lefteqn{ \sum_k p_k (I - aP)_{k, l-1}}\\&&
 = (1-a) p_{l-1} + aq_{l-1}p_{l-1} - aq_{l-2}p_{l-2}\alpha_{l-2,l-1} - a q_l p_l \alpha_{l,l-1} \\&&
= (1-a) \left\{p_{l-1} + \frac{(i_l - i_{l-2})}{(i_l - i_{l-1})(i_{l-1} - i_{l-2})}\sum_{k=1}^{l-2} (i_{l-1} - i_k) p_k \right. \\&&  \hspace{5mm} \left.  - \frac{1}{ i_{l-1} - i_{l-2} } \sum_{k=1}^{l-3} (i_{l-2} - i_k)p_k \right. \\&& \left. \hspace{20mm}  - \frac{1}{i_l - i_{l-1}}\left( e_0 - \sum_{k=1}^l i_k p_k - i_l + i_l\sum_{k=1}^l p_k \right) \right \} \\&& = (1-a)\left\{\frac{ i_l - e_0 }{ i_l - i_{l-1}} + \frac{i_l - i_l}{i_l - i_{l-1}}   p_l + \left(1 + \frac{(i_{l+1} - i_l)(i_{l-1} - i_l)}{(i_{l+1} - i_l)(i_l - i_{l-1})}\right ) p_{l-1}\right. \\&& \hspace{5mm} \left. + p_{l-2}\left( \frac{i_l - i_{l-2}}{i_l - i_{l-1}} - \frac{i_l - i_{l-2}}{i_l - i_{l-1}} \right) \right. \\&& \left. \hspace{10mm} + \sum_{k=1}^{l-3} p_k \left(\frac{(i_l - i_{l-2})(i_{l-1} - i_k)}{(i_l - i_{l-1})(i_{l-1} - i_{l-2})} - \frac{i_{l-2} - i_k}{(i_{l-1} - i_{l-2})} - \frac{i_l - i_{k}}{i_l - i_{l-1}} \right)\right \}\\&& = (1-a) \left\{ \frac{i_l - e_0}{i_l - i_{l-2}} \right\}
\end{eqnarray*}

\noindent as required. The computation giving 

\[ \sum_k p_k (I - aP)_{kl} = (1-a)\left\{ \frac{e_0 - i_{l-1}}{i_l - i_{l-1}}\right\} \]

\noindent is similar.  Lemma~\ref{lmmone} follows. \qed  
 
 \paragraph{Proof of theorem~\ref{thgt}} Let $\underline{q} = (q_1, \ldots q_{M}) \in
\{0\} \times (0,1]^{M-2} \times \{0\}$ and let $P(\underline{q})$ be the matrix  defined according to
definition~\ref{defdtgd}. Let $X$ denote the discrete time Markov chain with one step transition matrix given by $P(\underline{q})$ and transitions defined in definition~\ref{defdtgd} and let $\tau$ denote a random time, independent of $X$, with probability function 

\[ p(\tau = k) = a^k (1-a) \qquad k = 0,1,2, \ldots \]

\noindent Let $\mathbb{P}$ denote expectation with respect to both the random walk and the independent time $\tau$.  Then the matrix $P(\underline{q})$ provides a solution if and only if for each $j \in \{1, \ldots, M\}$ 

\begin{eqnarray}\label{eqefff}\lefteqn{   p_j =
\mathbb{P}(X_\tau = i_j | X_0 = e_0) }\\&& =  \nonumber
\frac{e_0 - i_{l-1}}{i_l - i_{l-1}}(1-a)\sum_{k=0}^\infty ((aP)^k)_{lj} +
\frac{i_l - e_0}{i_l - i_{l-1}}(1-a)\sum_{k=0}^\infty ((aP)^k)_{l-1,j} \qquad i_{l-1}
  < e_0 \leq i_l    
 \end{eqnarray}
 
\noindent where $(aP)^k = a^k P^k$ and, with $P^k$, multiplication is in the sense of matrix multiplication.  Set

\begin{equation}\label{eqgeepee} G(\underline{q}) = \sum_{k=0}^\infty (aP(\underline{q}))^k,
 \end{equation}

\noindent then $G(\underline{q})$ is well defined for $a \in (0,1)$ and satisfies

\begin{equation}\label{eqgeedef} G(I-aP) = (I - aP)G = I.\end{equation}

\noindent Therefore $\underline{q}$ provides a solution if and only if  $G(\underline{q})$ satisfies

\begin{equation}\label{eqmg2}  
p_j = \frac{e_0 - i_{l-1}}{i_l - i_{l-1}}(1-a)G_{lj}(\underline{q})+ \frac{i_l -
e_0}{i_l - i_{l-1}}(1-a)G_{l-1,j}(\underline{q}) \qquad j=1,\ldots, M  \qquad i_{l-1} < e_0 \leq i_l  \end{equation}

\noindent  It follows that  $\underline{q} \in \{0\} \times (0,1]^{M-2} \times \{0\}$  provides
a solution if and only if

\[ \sum_{k=1}^M p_k (I - aP(\underline{q}))_{k,j} = \left\{\begin{array}{ll} \frac{e_0
- i_{l-1}}{i_l - i_{l-1}}(1-a) & j = l\\ \frac{i_l - e_0}{i_l - i_{l-1}}(1-a) &
j = l-1 \\ 0 & j \neq l-1 \quad \mbox{or} \quad l  \end{array} \right. \] 

\noindent if $i_{l-1} < e_0 \leq i_l$, $i_{l-1}, i_l \in S$.  By
lemma~\ref{lmmone}, this system of equations has a unique solution, given by
equation (\ref{eqqq}). The    unique solution satisfies $\underline{q} \in \{0\} \times (0,1]^{M-2} \times \{0\}$ if and only if inequality (\ref{eqqbd}) holds. The proof is complete. \qed 

\section{The fixed point theorem and proof of theorem~\ref{thnb}}\label{secthnb}

In this section, theorem~\ref{thnb} is proved, establishing conditions under which there exists a DMRW (definition~\ref{defdtgd}) which, stopped at an independent negative binomial time, has a specified probability distribution. 
\paragraph{Description} Central to the proof is that a negative binomial variable $\tau \sim NB(r,a)$ is the sum of $r$ independent geometric variables, each with parameter $a$. If $\tau \sim NB(r,a)$, then
\begin{eqnarray*} \mathbb{P}(X_\tau = i_k | X_0 = i_j) &=& \sum_{n=0}^\infty \mathbb{P}(X_n = i_k | X_0 = i_j) \mathbb{P}(\tau = n)\\
&=& \sum_{n=0}^{\infty}(P^n(\underline{q}))_{jk} \left(\begin{array} {c}n+r-1\\n\end{array}\right)a^n (1-a)^r \\
&=& (1-a)^r \sum_{n=0}^\infty \left(\begin{array}{c}n+r-1\\n\end{array}\right)(aP^n(\underline{q}))_{jk} \\&=& 
 (1-a)^r (I - aP(\underline{q}))^{-r}_{jk} \\
 \end{eqnarray*}

\noindent where $P^n$ is taken in the sense of matrix multiplication of the one step transition matrices. It follows that if $\underline{v}$ is the distribution of $X_{0+}$ and $\underline{p}$ is the target distribution, then
\[ \underline{p} = (1-a)^r \underline{v} (I - aP(\underline{q}))^{-r}.\] If the process is started from point $e_0$ and $e_0 \not \in S$, definition~\ref{defdtgd} describes the rules for the first jump at time $0$; $\underline{v}$ gives the distribution at time $0+$.     The system $\underline{p}(I - aP(\underline{q}))^r = (1-a)^r \underline{v}$ is a system of degree $r$ polynomial equations; existence of solution is not straightfoward, let alone uniqueness, or any explicit expression for the solution. This article limits itself to existence. Existence of an $a_0 < 1$ such that for all $a \in [a_0,1)$ there exists a  suitable $\underline{q}$ is established as follows. The equation can be expressed as   $(1-a)^{-(r-1)}\underline{p} (I - aP(\underline{q}))^{r-1} = (1-a)  \underline{v} (I - aP(\underline{q}))^{-1}$.   In order to use theorem~\ref{thgt}, it has to be shown that $(1-a)^{-(r-1)}\sum_j (\underline{p} (I - aP(\underline{q}))^{r-1})_j = 1$ and $(1-a)^{-(r-1)} \sum_j i_j  (\underline{p} (I - aP(\underline{q}))^{r-1})_j = e_0(\underline{p})$. This is the subject of lemma~\ref{lmmaitchprex0}. Once this is established, it then follows, by theorem~\ref{thgt}, that any solution $\underline{q}$ satisfies

\[ \left\{\begin{array}{ll} q_j = \left(\frac{1}{a} - 1 \right) {\cal F}   \left((1-a)^{-(r-1)}\underline{p}(I - aP(\underline{q}))^{r-1} ,j \right) &  j \in \{2, \ldots, M-1\} \\ 0 & j = 1, M \end{array}\right. \]

\noindent  Existence of a range $[a_0,1)$ such that for any $a \in [a_0,1)$ there exists a solution $\underline{q}$  is now essentially a fixed point theorem. Let $\underline{\lambda} = \left(\frac{a}{1-a}\right)\underline{q}$. Then $(1-a)^{-1}(I - aP(\underline{q})) = N(\underline{\lambda})$ where the matrix $N$ is defined below. It follows that solutions satisfy $\underline{\lambda} = {\cal F}(\underline{p} N^{r-1} (\underline{\lambda}))$ and  theorem~\ref{thfpt},  gives existence of $\underline{\lambda} \in \{0\} \times {\bf R}^{M-2} \times \{0\}$ that satisfies this equation, together with existence of a positive lower bound and finite upper bound for the components of the fixed point. It follows that  $\underline{q} = \left(\frac{1}{a} - 1 \right) \underline{\lambda}$ provides a solution for $a \in [a_0,1)$ where $a_0$ satisfies $\left(\frac{1}{a_0} - 1 \right) \max_j \lambda_j = 1$, so that $\underline{q} \in \{0\} \times [0,1]^{M-2} \times \{0\}$.

Lemmas~\ref{lmmninv} and \ref{lmmnr1} are technical lemmas used in the proof of theorem~\ref{thfpt}. Lemma~\ref{lmmninv}, by finding a probabilistic expression for the entries of the matrix $N^{-(r-1)}$ (defined below) that appears in the proof of theorem~\ref{thfpt}, proves that the entries are non negative and that the rows sum to $1$. Lemma~\ref{lmmnr1} shows that for any $p \geq 1$, the rows of $N^p$ sum to $1$. Lemma~\ref{lmmnm120} deals with the columns of $N^{-p}$ for $p \geq 1$, showing that they tend to zero if the corresponding component of $\underline{\lambda}$ goes to infinity. 

Let $\underline{h}(\underline{p}, \underline{\lambda}) = \underline{p} N^{r-1} (\underline{\lambda})$, the notation used below. To construct a fixed point theorem, firstly the function $\underline{h}$ has to be modified to form a function $\underline{\tilde{h}}^{(\epsilon)}$ where the entries are bounded from above and below to ensure existence of a fixed point for the modified problem. The difficulty is with letting $\epsilon \rightarrow 0$ and showing that the components of $\underline{\lambda}^{(\epsilon)}$ are bounded from above and also bounded from below by a constant strictly greater than $0$, so that any limit point is a solution to the fixed point problem. Lemmas~\ref{lmmninv}, \ref{lmmnr1} and~\ref{lmmnm120} play an integral part in this, their roles are described more fully below. 

\paragraph{Notations} For $\underline{\lambda} \in \{0\} \times \mathbb{R}_+^{M-2} \times \{0\}$, define the $M \times M$ matrix 
$N (\underline{\lambda})$ by

\begin{equation}\label{eqndef} \left\{\begin{array}{ll}  N_{j,j} = 1 + \lambda_j
& j = 1,\ldots, M\\
N_{j,j+1} = - \lambda_j \alpha_{j,j+1}  &  j = 1,\ldots, M-1\\
 N_{j,j-1} = - \lambda_j \alpha_{j,j-1} &  j = 2,\ldots, M\\
N_{j,k} = 0 & \mbox{otherwise}.\end{array} \right. \end{equation}

\noindent where $\alpha_{j,k}$ is defined in equation (\ref{eqalpha}). For an $M$-row vector $\underline{b} \in {\bf R}^M$ and $\underline{\lambda} \in \{0\} \times \mathbb{R}_+^{M-2} \times \{0\}$, define $\underline{h} = (h_1, \ldots, h_M)$ by 

\begin{equation}\label{eqaitchemque}
 \underline{h}(\underline{b},\underline{\lambda}) = \underline{b} N^{r-1}(\underline{\lambda})
\end{equation}

\noindent Let $P(\underline{q})$ be defined as in equation (\ref{eqPdef}) and note that

\begin{equation}\label{eqnrep} \frac{1}{1-a}(I - aP(\underline{q})) =
N\left(\left(\frac{a}{1-a}\right)\underline{q} \right).
\end{equation}

\noindent This is easily seen;

\[ \frac{1}{1-a}(I - aP(\underline{q}))_{j,j} = \frac{1}{1-a}((1-a) + aq_j) = 1 +
\left(\frac{a}{1-a}\right)q_j \]
\[ \frac{1}{1-a}(I - aP(\underline{q}))_{j,j+1} =
-\alpha_{j,j+1} \left(\frac{a}{1-a}\right) q_j \]
\[ \frac{1}{1-a}(I-aP(\underline{q}))_{j,j-1} = - \alpha_{j,j-1}
\left(\frac{a}{1-a}\right)q_j.\]

\noindent It follows that 

\begin{equation}\label{eqha} h_j \left (\underline{p}, \left(\frac{a}{1-a}\right)\underline{q} \right
) = \frac{1}{(1-a)^{r-1}}\sum_{k=1}^M p_k ((I - aP(\underline{q}))^{r-1})_{k,j}.
\end{equation}

\noindent The following lemma shows that when $\underline{p}$ is a probability (that is $p_i \geq 0$ for each $i = 1,\ldots, M$ and $\sum_{j=1}^M p_j = 1$), then for each $\underline{\lambda}$, the components of $\underline{h}(\underline{p}, \underline{\lambda})$ sum to $1$ and
 $\sum_j i_j h_j(\underline{p},\underline{\lambda}) = \sum_j i_j p_j$. This means that if, furthermore, all the entries of $\underline{h}$ are positive, then it is a probability mass function. When taken over $S = \{i_1, \ldots, i_M\}$, it has  the same expectation as $\underline{p}$. 

\begin{Lmm}\label{lmmaitchprex0} Let $\underline{p}$ be an $M$ row vector that satisfies $p_i \geq 0$ for each $i \in \{1, \ldots, M\}$ and $\sum_{j=1}^M p_j = 1$.  Let $e_0(\underline{p}) = \sum_{j=1}^M i_j p_j$, let $\underline{\lambda} \in
\{0\} \times \mathbb{R}_+^{M-2} \times \{0\}$ and let $\underline{h}$ be defined in equation (\ref{eqaitchemque}). Then
\begin{equation}\label{eqss1} \sum_{j=1}^M h_j ( \underline{p}, \underline{\lambda}  ) = 1
\end{equation}
\noindent and
\begin{equation}\label{eqsx0} \sum_{j = 1} i_j h_j ( \underline{p} , \underline{\lambda} ) = e_0(\underline{p}).
\end{equation}
\end{Lmm}

\paragraph{Proof}  For equation (\ref{eqss1}), let $\underline{q} = \frac{\underline{\lambda}}{\max_{j
\in \{1,\ldots, M\}} \lambda_j}$, let $\lambda^* = \max_{j \in \{1, \ldots, M\}} \lambda_j$ and let $a = \frac{\lambda^*}{1 + \lambda^*}$, so that  $\frac{a}{1-a} \underline{q} =
\underline{\lambda}$. Using equation (\ref{eqha}),

\begin{eqnarray*} \sum_{j=1}^M h_j ( \underline{p}, \underline{\lambda} ) &=&  
\frac{1}{(1-a)^{r-1}}\sum_{k=1}^M p_k \sum_{i=0}^{r-1} (-1)^i
\left(\begin{array}{c} r - 1 \\ i \end{array}\right ) a^i \sum_{j=1}^M
(P^i)_{k,j} \\ &=& \frac{1}{(1-a)^{r-1}}\sum_{i=0}^{r-1} (-1)^i
\left(\begin{array}{c} r - 1 \\ i \end{array}\right ) a^i = 1.\end{eqnarray*}

\noindent For equation (\ref{eqsx0}), let $X_i$ denote the position at time $i$
of the Markov chain with one step transition matrix $P$ given in definition
(\ref{defdtgd}) determined by $\underline{q}$  and recall that $X$ is a martingale, so that 
\[ \mathbb{E}[X_i | X_0 = i_k] = i_k.\] 

\noindent Then, defining $a$ and $\underline{q}$ as above, 
 
\begin{eqnarray*} \sum_{j=1}^M i_j h_j (\underline{p}, \underline{\lambda} ) &=&
\frac{1}{(1-a)^{r-1}}\sum_{j=1}^M i_j \sum_{k=1}^M p_k \sum_{i=0}^{r-1} (-1)^i
\left(\begin{array}{c} r -1 \\ i \end{array}\right) a^i (P^i)_{ k, j} \\
&=&  \frac{1}{(1-a)^{r-1}}  \sum_{k=1}^M p_k \sum_{i=0}^{r-1} (-1)^i
\left(\begin{array}{c} r - 1 \\ i \end{array}\right) a^i \mathbb{E}[X_i | X_0 = i_k]\\
&=& \frac{1}{(1-a)^{r-1}}   \sum_{k=1}^M i_k p_k \sum_{i=0}^{r-1} (-1)^i
\left(\begin{array}{c} r - 1 \\ i \end{array}\right) a^i  \\
&=& \sum_{k=1}^M i_k p_k = e_0 (\underline{p}).
 \end{eqnarray*} \qed
 
 \noindent The following lemma gives a useful representation of $\underline{p}$ in terms of
$\underline{h}$ and gives useful properties of the inverse $N^{-1}$ that will be used later. The crucial properties are that elements of $N^{-1}$ are all non negative and each row sums to $1$, giving control on the entries. The characterisation of equation (\ref{eqninv}) will be used in the following way: if $a  < 1$, then $\mathbb{E}[\tau(a )] < +\infty$. If $q_j = 0$ for some $j$, then $P_{j,j} = 1$ and $\mathbb{P}( X_{\tau(a)} = i_{m_1} | X_0 = i_{m_2}) = 0$ for $m_1 < j \leq m_2$ or $m_1 > j \geq m_2$, which will be used to show that a contradiction is obtained if $\lambda_j^{(\epsilon)} \stackrel{\epsilon \rightarrow 0}{\longrightarrow} 0$, the contradiction being that the process cannot get past site $j$. The fixed point theorem~\ref{thfpt} starts by finding a fixed point $\underline{\lambda}^{(\epsilon)}$ to an approximate equation and then the limit $\epsilon \rightarrow 0$ is examined. It is shown, using equation (\ref{eqpchar}) and lemma~\ref{lmmninv} that $\min_j \lambda^{(0)}_j > 0$ for any limit point $\underline{\lambda}^{(0)}$.  

\begin{Lmm} \label{lmmninv} Let $r \geq 2$. 
  Let $\underline{p}= (p_1, \ldots, p_M)$ be a probability mass function satisfying $\inf_{j \in \{1,\ldots, M\}} p_j > 0$. Recall that $\underline{h}(\underline{p}, \underline{\lambda}) = \underline{p}N^{r-1}(\lambda)$ where $N$ is defined by equation
(\ref{eqndef}). Let $\underline{\lambda} \in \{0\} \times \mathbb{R}_+^{M-2} \times \{0\}$. Let
$\lambda^* = \sup_j \lambda_j$ and let $a_0 = \frac{\lambda^*}{1 + \lambda^*}$. For any $a \in [a_0,1)$, let $\tau(a)$ denote a random time with
probability function 
 
 \[ \mathbb{P} (\tau(a) = k) = (1-a)^{r-1} \left(\begin{array}{c} r -2 + k \\ k
\end{array}\right ) a^k \qquad k = 0,1,2, \ldots \]
 
 \noindent and let $X$ denote a DMRW (definition~\ref{defdtgd}), where the initial transition at time $0$ and subsequent one step transitions are defined in definition~\ref{defdtgd}, independent of $\tau(a)$,
with parameter vector $\underline{q} = \left(\frac{1}{a} - 1\right)\underline{\lambda}$. Then
 
 \begin{equation}\label{eqninv}
  (N^{-(r-1)})_{jk} = \mathbb{P}(X_{\tau(a)} = i_k | X_0 = i_j).
 \end{equation}
 
 \noindent It follows that $(N^{-(r-1)})_{j,k} > 0$ for all $\underline{\lambda} \in \{0\}
\times \mathbb{R}_+^{M-2} \times \{0\}$ (that is $\lambda_2, \ldots, \lambda_{M-1}$ all strictly positive) and that 
 
\begin{equation}\label{eqinvsum} \sum_k (N^{-1})_{jk} = 1 \qquad \forall  j \in
\{1, \ldots, M\}.
\end{equation}

\noindent It also follows that

 \begin{equation}\label{eqp} p_k = \sum_{j = 1}^M h_j (\underline{p}, \underline{\lambda} )
\mathbb{P}(X_{\tau(a)} = i_k | X_0 = i_j).\end{equation}
\end{Lmm}

\paragraph{Proof}  For such a choice of $a$ and $\underline{q}$, it follows from equation
(\ref{eqnrep}) that

\[ N^{-(r-1)}(\lambda) = (1-a)^{r-1} (I - aP(\underline{q}))^{-(r-1)} = (1-a)^{r-1}\sum_{k =
0}^\infty \left(\begin{array}{c} r - 2 + k \\ k \end{array}\right ) (aP(\underline{q}))^k\]

\noindent so that

\begin{equation}\label{eqninvexp} N^{-(r-1)}_{m_1,m_2}(\underline{\lambda}) =
(1-a)^{r-1} \sum_{k=0}^\infty \left(\begin{array}{c} r - 2 + k\\ k
\end{array}\right) a^k \mathbb{P}(X_k = i_{m_2} | X_0 = i_{m_1})\end{equation}

\noindent and equation (\ref{eqninv}) follows from the definition of $\tau(a)$;
equation (\ref{eqinvsum}) follows from summing over $m_2$ in equation
(\ref{eqninvexp}).
It follows from the definition of $\underline{h}$ given by equation (\ref{eqaitchemque}) that 
\[ \underline{p} = \underline{h}(\underline{p}, \underline{\lambda}) N^{-(r-1)}(\underline{\lambda})\] and hence that, for $k \in \{1, \ldots, M\}$, 

\begin{eqnarray*} p_k &=& (1-a)^{r-1} \sum_{i=0}^\infty
\left(\begin{array}{c} r + i - 2 \\ i \end{array}\right) a^i \sum_{j=1}^M
h_j (\underline{p}, \underline{\lambda} ) (P^i)_{j,k} \\
&=&  (1-a)^{r-1} \sum_{i=0}^\infty \left(\begin{array}{c} r + i - 2 \\ i
\end{array}\right) a^i \sum_{j=1}^M h_j (\underline{p}, \underline{\lambda} ) \mathbb{P}(X_i = i_k | X_0 = i_j)
\\
&=& \sum_{j = 1}^M h_j (\underline{p},\underline{\lambda} ) \mathbb{P}(X_{\tau(a)} = i_k | X_0 = i_j) 
\end{eqnarray*}

\noindent thus establishing equation (\ref{eqp}) and completing the proof of
lemma~\ref{lmmninv}. \qed 

\begin{Lmm}\label{lmmnr1} For all $r \geq 1$, $\sum_j (N^r)_{k,j} = 1$ for each
$k$.
 \end{Lmm}
 
 \paragraph{Proof} By construction, this is clear for $r = 1$;
 \[ \sum_j N_{k,j} = \left\{ \begin{array}{ll} -\alpha_{k,k-1} \lambda_k + (1 +
\lambda_k) - \alpha_{k,k+1}\lambda_k = 1 & k \in \{2,\ldots, M-1\} \\
  1 & k \in \{1,M\} \end{array}\right. 
 \]
For $r \geq 2$, assume that the result is true for $r - 1$, then
\[ \sum_j (N^r)_{k,j}  = \sum_{i,j} (N^{r-1})_{k,i} N_{ij} = \sum_i (N^{r-1})_{k,i}  =
1.\]
\noindent  The result follows by induction. \qed \vspace{5mm}

\noindent Lemma~\ref{lmmnm120} gives important limiting behaviour of the entries of $N^{-p}$ for each $p$. Since it has been established by lemma~\ref{lmmninv} that the entries are non negative and the rows sum to $1$, lemma~\ref{lmmnm120} is a key step to showing that $\max_j \sup_\epsilon \lambda^{(\epsilon)}_j$ is bounded, where $\underline{\lambda}^{(\epsilon)}$ is a fixed point for the approximating problem in theorem~\ref{thfpt}. 

\begin{Lmm}\label{lmmnm120}
 If $\lambda_j \rightarrow +\infty$ then $(N^{-1})_{.j} \rightarrow \underline{0}$ and
consequently $(N^{-r})_{.j} \rightarrow \underline{0}$ (in the limit, each entry of the $j$th column is identically equal to zero) for any integer $r \geq 1$.
\end{Lmm}

\paragraph{Proof of Lemma~\ref{lmmnm120}}
Let $\beta_{kj} = (N^{-1})_{kj}$. Then $0 \leq \beta_{kj} \leq 1$ by
lemma~\ref{lmmninv} and $\beta$ satisfies the following system of equations:

\[ -\alpha_{k,k-1} \beta_{k-1,j} + \left (1 + \frac{1}{\lambda_k} \right )\beta_{k,j} -
\alpha_{k,k+1}\beta_{k+1,j} = \left\{\begin{array}{ll} 0 & k \neq j \\
\frac{1}{\lambda_j} & k = j. \end{array}\right. \]

\noindent where $\alpha_{k,k-1}\beta_{k-1,j} = 0$ when $k = 1$ by definition and
$\alpha_{k,k+1}\beta_{k+1,j} = 0$ when $k = M$ by definition. \vspace{5mm}

\noindent From the a priori bounds on
$\beta_{k,j}$ from the previous lemma ($0 \leq \beta_{k,j} \leq 1$), it follows  that
$\beta_{k,j} \stackrel{\lambda_j \rightarrow +\infty}{\longrightarrow} 0$ for
all $k = 1,\ldots, M$. \qed\vspace{5mm}

\noindent Having stated and proved the preparatory lemmas, the fixed point theorem may be stated and proved.

\begin{Th}[Fixed Point Theorem]\label{thfpt}
Let $\underline{p} = (p_1, \ldots, p_M)$ be a probability mass function satisfying $\min_{k \in \{1, \ldots, M\}} p_k > 0$ over a set of points $S = \{i_1, \ldots, i_M\}$, with expectation 
 $\sum_{j=1}^M i_j p_j = e_0(\underline{p}) \in \mathbb{R}$. Recall that 
$\underline{h}(\underline{p},\underline{\lambda})$ is defined as in equation (\ref{eqaitchemque}) and let ${\cal F}$ be defined as in equation (\ref{eqefdef}). Assume that $\underline{p}$ is fixed. There exists a
$\underline{\lambda} \in \{0\} \times \mathbb{R}_+^{M-2} \times \{0\}$ satisfying

\begin{equation}\label{eqfpeqn} \underline{\lambda} =  {\cal F}\left (\underline{h}\left ( \underline{p}, \underline{\lambda} \right )
\right ).
\end{equation}

\noindent  This point satisfies

\begin{equation}\label{eqlambdconc}  0 < \min_{j \in \{2,\ldots,M-1\}} \lambda_j
\leq \max_{j \in \{2, \ldots, M-1\}} \lambda_j < +\infty\end{equation}

\noindent and

\begin{equation}\label{eqhconc} \min_{j \in \{1, \ldots, M\}} h_j (\underline{p},\underline{\lambda} ) >
0.
\end{equation}

\end{Th}

\paragraph{Proof of theorem~\ref{thfpt}. Part 1: Fixed point for an approximating problem}

Let $\underline{h}$ satisfy equation (\ref{eqaitchemque}). For $0 < \epsilon < 1$, let
\begin{equation}\label{eqcle} C(\underline{\lambda}, \epsilon) = \sum_{j=1}^M
\left(\frac{h_j (\underline{p},\underline{\lambda})}{\sum_{k=1}^M (h_k (\underline{p},\underline{\lambda} ) \vee \epsilon)}
\vee \epsilon\right)\end{equation}

\noindent and let $\underline{\tilde{h}}^{(\epsilon)}$ be defined as 

\begin{equation}\label{eqaitcheps} \tilde{h}^{(\epsilon)}_j ( \underline{p} , \underline{\lambda} ) =
\frac{1}{C(\underline{\lambda},\epsilon)}\left(\frac{h_j (\underline{p}, \underline{\lambda} )}{\sum_{k = 1}^M (h_k 
( \underline{p}, \underline{\lambda} ) \vee \epsilon)} \vee \epsilon \right), \end{equation}

\noindent so that $\sum_j \tilde{h}_j ^{(\epsilon)}(\underline{p}, \underline{\lambda} ) = 1$.   From
equation (\ref{eqss1}) lemma~\ref{lmmaitchprex0}, it follows directly that 
\[\sum_k (h_k (\underline{p}, \underline{\lambda} ) \vee \epsilon ) \geq 1,\]

\noindent since for each $k$, $h_k (\underline{p}, \underline{\lambda} ) \vee \epsilon \geq h_k (\underline{p}, \underline{\lambda} )$ and $\sum_j h_k (\underline{p}, \underline{\lambda} ) = 1$. \vspace{5mm}

\noindent From this, it follows directly that $C(\underline{\lambda},\epsilon) \geq 1$, which follows because $\epsilon \geq \frac{\epsilon}{\sum_k (h_k (\underline{p}, \underline{\lambda} ) \vee \epsilon )}$ and hence 

\begin{eqnarray*} C(\lambda, \epsilon) &=& \sum_{j=1}^M
\left(\frac{h_j (p,\lambda )}{\sum_{k=1}^M (h_k (\underline{p},\underline{\lambda} ) \vee \epsilon)}
\vee \epsilon\right) \geq  \sum_{j=1}^M \left (
 \frac{h_j (p,\lambda )}{\sum_{k=1}^M (h_k (\underline{p},\underline{\lambda} ) \vee \epsilon)}  
\vee  \frac{\epsilon}{\sum_{k=1}^M (h_k (\underline{p},\underline{\lambda} ) \vee \epsilon) }  \right ) \\ &=& \frac{1}{\sum_{k=1}^M (h_k (\underline{p},\underline{\lambda} ) \vee \epsilon)} \sum_{k=1}^M (h_k (\underline{p},\underline{\lambda} ) \vee \epsilon) = 1.
\end{eqnarray*}

\noindent Furthermore, for each $j \in \{1, \ldots, M\}$, 

\[ \frac{h_j(\underline{p}, \underline{\lambda})}{\sum_{k=1}^M (h_k (\underline{p}, \underline{\lambda}) \vee \epsilon)} \leq 1\] so that, for $0 < \epsilon < 1$, $\sum_{j=1}^M
\left(\frac{h_j (p,\lambda )}{\sum_{k=1}^M (h_k (\underline{p},\underline{\lambda} ) \vee \epsilon)}
\vee \epsilon\right) \leq \sum_{j=1}^M (1 \vee \epsilon) = M$. 
It therefore follows that 
\[ 1 \leq C(\underline{\lambda},\epsilon) = \sum_{j=1}^M \left(\frac{h_j (\underline{p}, \underline{\lambda})}{\sum_{k = 1}^M (h_k (\underline{p},\underline{\lambda} ) \vee \epsilon)} \vee \epsilon
\right) \leq M.\]  

\noindent Since $\frac{h_j (\underline{p}, \underline{\lambda} )}{\sum_{k = 1}^M (h_k 
( \underline{p}, \underline{\lambda} ) \vee \epsilon)} \vee \epsilon \geq \epsilon$,  it follows that 

\begin{equation}\label{eqhlb}
 \inf_\lambda \min_j \tilde{h}^{(\epsilon)}_j(\underline{p}, \underline{\lambda}) \geq
\frac{\epsilon}{M}.
\end{equation}

\noindent Let

\[ {\cal A}^{(\epsilon)} (\underline{\lambda}, \underline{p})(j) :=  {\cal F} \left (
\underline{\tilde{h}}^{(\epsilon)}\left ( \underline{p}, \underline{\lambda} \right ), j \right ) \qquad j \in \{1,
\ldots, M\}\]

\noindent  where ${\cal F}$ is defined in equation (\ref{eqefdef}).  \vspace{5mm}

\noindent Recall the definition of ${\cal L}$  in equation (\ref{eqeldef}). It
follows directly that for any probability distribution $\underline{p}$ over $\{1, \ldots, M\}$,   

\begin{equation}\label{eqlmax} {\cal L}(\underline{p},j) \leq \max_{j \in \{2,\ldots,
M-1\}}\frac{(i_{j+1} - i_{j-1})(i_M - i_1)}{(i_{j+1} - i_j)(i_j - i_{j-1})} <
+\infty,\end{equation}

\noindent and hence that 

\[ 0 \leq {\cal A}^{(\epsilon)}(\underline{\lambda} , \underline{p})(j) \leq \frac{M}{\epsilon} \max_{j
\in \{2,\ldots, M-1\}}\frac{(i_{j+1} - i_{j-1})(i_M - i_1)}{(i_{j+1} - i_j)(i_j
- i_{j-1})} =: C^{(\epsilon)} < +\infty.\]

\noindent for a constant $C^{(\epsilon)}$. Then, for fixed $\underline{p}$, ${\cal A}^{(\epsilon)}( .,\underline{p}) : \{0\} \times {\bf
R}_+^{M-2} \times \{0\} \rightarrow \{0\} \times [0,C^{(\epsilon)}]^{M-2} \times
\{0\}$. The space $\{0\} \times [0, C^{(\epsilon)}]^{M-2}\times \{0\}$ is a
closed convex subset of $\mathbb{R}^M$. The following consideration shows that, for
fixed $\epsilon > 0$,  the map ${\cal A}^{(\epsilon)}(.,\underline{p})$ is continuous. From
the definition of $\underline{h}$ (equation (\ref{eqaitchemque})), it follows that  $\underline{h}(\underline{p},.)$
is a smooth function in $\lambda$. It therefore follows directly that
$\underline{\tilde{h}}^{(\epsilon)}(\underline{p},.)$ is continuous in $\underline{\lambda}$ and hence, using the lower
bound (\ref{eqhlb}), it follows directly from equation (\ref{eqefdef}) that for
any fixed $\epsilon > 0$, ${\cal A}^{(\epsilon)}(.,\underline{p})$ is  continuous in
$\underline{\lambda}$.  It follows that the mapping ${\cal A}^{(\epsilon)}(.,\underline{p}) : \{0\}
\times [0,C^{(\epsilon)}]^{M-2} \times \{0\} \mapsto  \{0\} \times
[0,C^{(\epsilon)}]^{M-2} \times \{0\}$  has a fixed point $\underline{\lambda}^{(\epsilon)}$
by Schauder's fixed point theorem, which states that if $K$ is a convex subset
of a topological vector space $V$ and $T$ is a continuous mapping of $K$ into
itself such that $T(K)$ is contained in a compact subset of $K$, then $T$ has a
fixed point.  

\paragraph{Part 2: showing $\sup_\epsilon \max_j \lambda_j^{(\epsilon)} < +\infty$.} 

Let $\underline{h}_\epsilon := \underline{\tilde{h}}^{(\epsilon)}(\underline{p}, \underline{\lambda}^{(\epsilon)})$. Then

\begin{equation}\label{eqlameps} \underline{\lambda}^{(\epsilon)} = {\cal F}(\underline{h}_{\epsilon})
 \end{equation}
 
 \noindent so that, from the definition of ${\cal L}$ (equation (\ref{eqeldef})),  

\begin{equation}\label{eqhlaml} h_{\epsilon,j} \lambda^{(\epsilon)}_j =  {\cal
L} \left (\underline{h}_\epsilon;j \right ) \qquad j = 1, \ldots, M.
 \end{equation}
 
 \noindent It follows from equation (\ref{eqlmax}) that 
 
 \begin{equation}\label{eqhlm} h_{\epsilon,j} \leq
\frac{1}{\lambda_j^{(\epsilon)}}\max_{j \in \{2,\ldots, M-1} \frac{(i_{j+1} -
i_{j-1})(i_M - i_1)}{(i_{j+1} - i_j)(i_j - i_{j-1})}
\stackrel{\lambda_j^{(\epsilon)} \rightarrow +\infty}{\longrightarrow}
0.\end{equation}

\noindent The aim of this part is to show firstly that $\sup_\epsilon \max_j
\lambda^{(\epsilon)}_j < +\infty$ and the next part that $\inf_\epsilon \min_j h_{\epsilon, j} > 0$. From this, it follows directly that any limit point $\lambda$ of
$\lambda^{(\epsilon)}$ satisfies equation (\ref{eqfpeqn}).  \vspace{5mm}

\noindent Set

\begin{equation}\label{eqph} \underline{p}^{(\epsilon)} = \underline{h}_\epsilon
N^{-(r-1)}(\underline{\lambda}^{(\epsilon)}),
\end{equation}

\noindent so that $\underline{h}_\epsilon = \underline{h}(\underline{p}^{(\epsilon)}, \lambda^{(\epsilon)})$. By
construction, since each $h_{\epsilon, j} > 0$ and $\sum_j h_{\epsilon, j} = 1$,
it follows using lemma~\ref{lmmninv}, where it is proved that $(N^{-(r-1)})_{ij}
\geq 0$ and $\sum_j (N^{-(r-1)})_{ij} = 1$,  that for each $j$,

\[ 0 \leq p^{(\epsilon)}_j \leq 1.\]

\noindent Furthermore,

\begin{equation}\label{eqspe} \sum_j p^{(\epsilon)}_j = \sum_{k}
h_{\epsilon,k}\sum_j (N^{-(r-1)})_{kj} = \sum_k h_{\epsilon, k} =
1.\end{equation}

\noindent Also, it follows from lemma~\ref{lmmnm120} that if
$\lambda_j^{(\epsilon)} \rightarrow +\infty$, then
$(N^{-(r-1)}(\lambda^{(\epsilon)}))_{.j} \rightarrow \underline{0}$. \vspace{5mm}

\noindent Using the bounds on the components of $\underline{h}_\epsilon$ together with this observation, it
follows from equation (\ref{eqph}) that if $\lambda_j^{(\epsilon)} \rightarrow
+\infty$, then $p^{(\epsilon)}_j \rightarrow 0$. \vspace{5mm}

\noindent From equation (\ref{eqss1}) lemma~\ref{lmmaitchprex0}, it follows
that 

\begin{equation}\label{eqsbd} K(\underline{\lambda},\epsilon) := \sum_{j} \left(h_j (\underline{p}, \underline{\lambda} ) \vee \epsilon \right) \geq 1.\end{equation}

\noindent  From the definition of $\underline{h}$ (equation (\ref{eqaitchemque})), it
follows that $\underline{h}(\underline{p}, \underline{\lambda})$ is linear in $\underline{p}$. Let

\begin{equation} \label{eqkepsdef} 
  C_\epsilon = C(\underline{\lambda}^{(\epsilon)},\epsilon) \qquad \mbox{and} \qquad K_\epsilon =
K(\underline{\lambda}^{(\epsilon)},\epsilon).
\end{equation}

\noindent Then, from the definition of $\underline{h}_\epsilon$,

\[ h_{\epsilon,j} =
\frac{1}{C_\epsilon K_\epsilon}
h_j (\underline{p}, \underline{\lambda}^{(\epsilon)} ) \vee \frac{\epsilon}{C_\epsilon}    = h_j \left(\frac{\underline{p}}{C_\epsilon K_\epsilon},
\underline{\lambda}^{(\epsilon)} \right) \vee \frac{\epsilon}{C_\epsilon}.\]

\noindent Note, from above, that $1 \leq C_\epsilon \leq M$ and $K_\epsilon \geq
1$. Let $N^{(\epsilon)}_{r-1}$ denote $N^{r-1}(\lambda^{(\epsilon)})$ with
column $k$ replaced by a column where each entry is $\epsilon K_\epsilon$ for each $k$ such that 
$h_{\epsilon,k} = \frac{\epsilon}{C_\epsilon}$. Since $\sum_j p_j = 1$, it
follows that 

\begin{equation}\label{eqhepsform} \underline{h}_\epsilon = \frac{1}{C_\epsilon K_\epsilon} \underline{p} N^{(\epsilon)}_{r-1}.
\end{equation}

\noindent Let $S_\epsilon = \{ \beta | h_{\epsilon, \beta} =
\frac{\epsilon}{C_\epsilon}\}$. Now, by construction, note that 

\begin{equation}\label{eqnmrm1rm1} (N^{-(r-1)}(\underline{\lambda}^{(\epsilon)})N^{(\epsilon)}_{r-1})_{m_1, m_2} =
\left\{ \begin{array}{ll} I(m_1, m_2) & m_2 \not \in S_\epsilon \\
\epsilon K_\epsilon & m_2 \in S_\epsilon \end{array} \right. \end{equation}

\noindent where $I(m_1 , m_2) = 1$ if $m_1 = m_2$ and $I(m_1, m_2)
= 0$ if $m_1 \neq m_2$. This follows because $\sum_j N^{-(r-1)}_{m_1, j}
= 1$ for each $m_1$. Let 
\begin{equation}\label{eqdeffeps} F_\epsilon = N^{-(r-1)}(\underline{\lambda}^{(\epsilon)})N^{(\epsilon)}_{r-1}.
\end{equation}

\noindent Directly from equation (\ref{eqdeffeps}), it follows that $N^{(\epsilon)}_{r-1} =
N^{(r-1)}(\underline{\lambda}^{(\epsilon)})F_\epsilon$ so that, from equation (\ref{eqhepsform}),  

\[ \underline{h}_\epsilon = \frac{1}{K_\epsilon C_\epsilon}\underline{p}
N^{(r-1)}(\lambda^{(\epsilon)})F_\epsilon \] 

\noindent and hence, from equation (\ref{eqph}),  

\begin{equation}\label{eqpeps} \underline{p}^{(\epsilon)} = \frac{1}{K_\epsilon C_\epsilon} \underline{p}
N^{(r-1)}(\underline{\lambda}^{(\epsilon)})F_\epsilon N^{-(r-1)}(\underline{\lambda}^{(\epsilon)}). 
\end{equation}

\noindent Let $\Lambda_\epsilon$ be the matrix such that
\[ \Lambda_{\epsilon ;m_1,m_2} = \left\{ \begin{array}{ll} 1 & m_1 = m_2, \qquad m_2 \not \in
S_\epsilon \\ 0 & \mbox{otherwise}. \end{array} \right. \]  

\noindent That is,
$\Lambda_\epsilon$ is a diagonal matrix with entries $\Lambda_{\epsilon;
m,m} = 1$ if $m \not \in S_\epsilon$ and $\Lambda_{\epsilon;
m,m} = 0$ if $m   \in S_\epsilon$, and let ${\cal J}_\epsilon$ denote the matrix such that 

\[ {\cal
J}_{\epsilon;m_1,m_2} = \left\{ \begin{array}{ll} 1 & m_2 \in S_\epsilon \\  0 & \mbox{otherwise}. \end{array}\right. \]
That is, ${\cal J}_\epsilon$ has columns of $1$s corresponding to elements of
$S_\epsilon$ and the remaining columns are columns of $0$s. Then, from equation (\ref{eqnmrm1rm1}) and the definition of $F_\epsilon$ in equation (\ref{eqdeffeps}), it follows that 

\[ F_\epsilon = \Lambda_\epsilon + \epsilon K_\epsilon {\cal J}_\epsilon.\]

\noindent From lemma~\ref{lmmnr1}, it follows that $\sum_j N^{(r-1)}_{m , j}
= 1$ for each $m$, from which it follows that

\[N^{(r-1)}(\underline{\lambda}^{(\epsilon)}) F_\epsilon N^{-(r-1)}(\underline{\lambda}^{(\epsilon)}) =
N^{(r-1)}(\underline{\lambda}^{(\epsilon)}) \Lambda_\epsilon
N^{-(r-1)}(\underline{\lambda}^{(\epsilon)}) + \epsilon K_\epsilon {\cal J}_\epsilon
N^{-(r-1)}(\underline{\lambda}^{(\epsilon)}) \]

\noindent and hence from equation (\ref{eqpeps}) that

\begin{equation}\label{eqpeps2} \underline{p}^{(\epsilon)} = \frac{1}{K_\epsilon C_\epsilon} \underline{p}
N^{(r-1)}(\underline{\lambda}^{(\epsilon)})\Lambda_\epsilon N^{-(r-1)}(\underline{\lambda}^{(\epsilon)}) + \frac{\epsilon}{C_\epsilon} \underline{p} {\cal J}_{\epsilon} N^{-(r-1)}(\underline{\lambda}^{(\epsilon)}) . 
\end{equation}

\noindent Since the elements of $N^{-(r-1)}$ are non negative and bounded by
$1$ and $C_\epsilon \geq 1$, it follows that 

\[  
\lim_{\epsilon \rightarrow 0} \frac{\epsilon }{ C_\epsilon} \underline{p}   {\cal J}_\epsilon
N^{-(r-1)}(\lambda^{(\epsilon)}) = 0.\]
 
\noindent For any invertible matrix ${\cal
M}$, the eigenvalues of ${\cal M}^{-1} {\cal A} {\cal M}$ are the same as the
eigenvalues of ${\cal A}$. It follows that the eigenvalues of
$N^{(r-1)}(\lambda^{(\epsilon)}) \Lambda_\epsilon
N^{-(r-1)}(\lambda^{(\epsilon)})$ are the eigenvalues of $\Lambda_\epsilon$; $0$
with multiplicity equal to the number of elements of $S_\epsilon$ and the
remaining eigenvalues all $1$.
\vspace{5mm}

\noindent It now follows directly that if $K_\epsilon \stackrel{\epsilon \rightarrow 0}{\longrightarrow} +\infty$, then 
$\underline{p}^{(\epsilon)} \stackrel{\epsilon \rightarrow 0}{\longrightarrow} \underline{0}$, contradicting the fact that $p^{(\epsilon)}_j \geq 0$ for each $j$ and $\sum_j
p^{(\epsilon)}_j = 1$ for all $\epsilon \in (0,1)$. \vspace{5mm}

\noindent It follows that 
\begin{equation} \label{eqske}\sup_\epsilon K_\epsilon < +\infty.
\end{equation}
 \vspace{5mm}

\noindent It follows from the definition of $K_\epsilon$ (equation (\ref{eqsbd}) and (\ref{eqkepsdef})) and the conclusion (\ref{eqske}) that 
\[\sum_{j=1}^M
h_j (\underline{p}, \underline{\lambda}^{(\epsilon)} ) \vee 0 \leq \sup_\epsilon K_\epsilon < +\infty.\] 

\noindent Since $\sum_j h_j(\underline{p}, \underline{\lambda}) = 1$
(lemma~\ref{lmmaitchprex0} equation (\ref{eqss1})) for all $\underline{\lambda} \in \{0\}
\times \mathbb{R}_+^{M-2} \times \{0\}$, it follows directly that  $\sup_\epsilon \max_j |h_j ( \underline{p}, \underline{\lambda}^{(\epsilon)} )| <
+\infty$.\vspace{5mm}

\noindent Let $ \lambda^{*(\epsilon)} = \max_j \lambda^{(\epsilon)}_j$ and let
$N^{*(\epsilon)} = \frac{1}{\lambda^{*(\epsilon)}}N(\underline{\lambda}^{(\epsilon)})$ (that
is, divide every element by $\lambda^{*(\epsilon)}$). Let $g_j =
\frac{1}{\lambda^{*(\epsilon)(r-1)}} h_j (\underline{p}, \underline{\lambda}^{(\epsilon)} )$. Then if
there is a sequence $\epsilon_n \rightarrow 0$ such that 
$\lambda^{*(\epsilon_n)} \stackrel{n \rightarrow +\infty}{\longrightarrow}
+\infty$, it follows that $g \stackrel{n \rightarrow +\infty}{\longrightarrow}
0$ and hence that, for any limit point $N^*$ of $N^{*(\epsilon_n)}$,

\[ 0 = \underline{p}N^{*(r-1)}.\]

\noindent But it follows from the construction of $N^*$ that the rank $\rho$ of
$N^*$ is the number of components of $\underline{\lambda}$ such that $\lim_{n \rightarrow
+\infty} \frac{\lambda_j^{(\epsilon_n)}}{\lambda^{*(\epsilon_n)}} > 0$, where
$\underline{\lambda}^{(\epsilon_n)}$ is a sequence that gives the limit point. This is
seen as follows: consider the lowest index $k_1$ such that $\lim_{n \rightarrow
+\infty}\frac{\lambda_{k_1}^{(\epsilon_n)}}{\lambda^{*(\epsilon_n)}} > 0$, then
in $N^*$, the limit, column $k_1 - 1$ will have exactly one entry; element
$N_{k_1,k_1-1}^*$ will be the only non-zero element of column $k_1$. Suppose
$k_1 < \ldots < k_\rho$ are the relevant indices, then the columns
$(N_{.,k_1-1}^*, \ldots, N_{., k_\rho - 1}^*)$ provide an upper triangular
matrix, with elements $N_{k_j, k_j - 1}^* \neq 0$ and $N_{p,k_j - 1} = 0$ for
all $p \geq k_j + 1$, proving that $N^*$ is of rank $\rho$. 

It follows that $N^{*(r-1)}$ is of rank $\rho$ and the non-zero rows of
$N^{*(r-1)}$ are those corresponding to the indices $k : \lim_{n \rightarrow
+\infty} \frac{\lambda_k^{(\epsilon_n)}}{\lambda^{*(\epsilon_n)}} > 0$. Since
the space spanned by the $\rho$ rows is of rank $\rho$, it follows that $p_k
= 0$ for each of these $p_k$, which is a contradiction (since, by hypothesis,
$p_k > 0$ for each $k$). \vspace{5mm}

\noindent Hence
\begin{equation}\label{eqsuplam} \sup_\epsilon \lambda^{*(\epsilon)} < +\infty. \end{equation} 

\paragraph{Part 3: showing $\inf_\epsilon \min_j \lambda^{(\epsilon)}_j > 0$.} Now suppose that $\lambda^{(\epsilon)}_j \rightarrow 0$ for some $j
\in \{2, \ldots, M-1\}$. Let $q_{\epsilon,j} =
\frac{\lambda^{(\epsilon)}_j}{\lambda^{*(\epsilon)}}$ and $a_\epsilon =
\frac{\lambda^{*(\epsilon)}}{1 + \lambda^{*(\epsilon)}}$ so that 
\[  \underline{q}_\epsilon = \left(\frac{1}{a_\epsilon} - 1 \right) \underline{\lambda}^{(\epsilon)}. \]

\noindent It follows from equation (\ref{eqsuplam}) that  $\sup_\epsilon a_\epsilon < 1$, where the inequality is strict. It follows from lemma~\ref{lmmninv} equation (\ref{eqninv}), because   $q_{\epsilon,j}
\rightarrow 0$ while $\sup_\epsilon a_\epsilon < 1$ (the inequality is strict),
that $N^{-(r-1)}_{kp}(\underline{\lambda}^{(\epsilon)}) \rightarrow 0$ for all $(k,p)$ such
that $k \leq j < p$ or $k \geq j > p$. This follows from the following consideration: let $\tau(a_\epsilon)$ denote a random time with probability function 

\[
\mathbb{P}(\tau(a_\epsilon) = k) = (1-a_\epsilon )^{r-1}\left(\begin{array}{c} r +
k - 2 \\ k \end{array}\right) a_\epsilon^k \qquad k = 0,1,2,
\ldots \] 

\noindent independent of the DMRW (definition~\ref{defdtgd}) with parameters  $\underline{q}_\epsilon$. Note that \[ \sup_\epsilon \mathbb{E}[\tau(a_\epsilon)] = \frac{(r-1)a^*}{1-a^*} <
+ \infty,\] where $a^* = \sup_\epsilon a_\epsilon$. If $q_{\epsilon,j} \rightarrow 0$
for some $m_1 \leq j < m_2$ where $m_1 < m_2$, then

\begin{eqnarray*} \lefteqn{\mathbb{P}(X_{\tau(a_\epsilon)}^{(\epsilon)} = i_{m_2} | X_{0}^{(\epsilon)} = i_{m_1}) = \sum_{k=0}^\infty \mathbb{P}(X_k^{(\epsilon)} = i_{m_2} | X_{0}^{(\epsilon)} = i_{m_1} , \tau(a_\epsilon) = k ) \mathbb{P}(\tau(a_\epsilon) = k)} \\ &&= 
\sum_{k=0}^\infty \sum_{k_1 = 0}^k \mathbb{P}(X_{k_1}^{(\epsilon)} = i_j|X_0^{(\epsilon)} = i_{m_1})\mathbb{P}(X_{k-k_1}^{(\epsilon)} = i_{m_2} |X_0^{(\epsilon)} = i_j) \mathbb{P}(\tau(a_\epsilon) = k).
\end{eqnarray*}

\noindent If $q_{\epsilon_n,j} \rightarrow 0$, then $\mathbb{P}(X_{k}^{(\epsilon_n)} = i_{m_2} |X_0^{(\epsilon_n)} = i_j) \rightarrow 0$ for all $k$. It follows from lemma~\ref{lmmninv} equation (\ref{eqninv}) that  
 \[  \mathbb{P}(X_{\tau(a_{\epsilon_n})}^{(\epsilon_n)} = i_{m_2} | X_{0}^{(\epsilon_n)} = i_{m_1} ) = (N^{-(r-1)})_{m_1, m_2}(\lambda^{(\epsilon_n)})
\rightarrow 0  \]

\noindent for all $(m_1, m_2)$ such  $m_1 \leq j <
m_2$. The proof that the result holds for $m_1 \geq j > m_2$ is similar.   

Furthermore, it follows from equation (\ref{eqhlaml}) (${\cal L}$ is defined in
equation (\ref{eqeldef})) that for any sequence with limit point $\underline{\lambda}^{(0)}$ such that 
$\lambda^{(\epsilon_n)}_j \rightarrow 0$ for some $j$, if
$j \leq l-1$ then $h_k (\underline{p}, \underline{\lambda}^{(0)} ) \leq 0$ for all $1 \leq k \leq j-1$
and if $j \geq l$ then $h_k (\underline{p}, \underline{\lambda}^{(0)} ) \leq 0$ for all $j+1 \leq k \leq
M$. \vspace{5mm}

\noindent It follows, using 
\begin{equation}\label{eqpchar} \underline{p} = \underline{h}(\underline{p}, \underline{\lambda}) N^{-(r-1)}(\lambda),\end{equation}

\noindent and taking $\underline{\lambda} = \underline{\lambda}^{(0)}$ and considering the zeroes of
$N^{-(r-1)}(\underline{\lambda}^{(0)})$ that if $j \leq l-1$, then $p_1 \leq 0, \ldots,
p_{j-1} \leq 0$, which is a contradiction. If $j \geq l$, then $p_{j+1}
\leq 0, \ldots, p_M \leq 0$, which is a contradiction.\vspace{5mm}

\noindent It follows that any limit point $\underline{\lambda}$ satisfies the bounds of
(\ref{eqlambdconc})  and consequently that 
\[ h_{0,1} > 0, \ldots h_{0,M} > 0,\] consequently that $\underline{h}( \underline{p} , \underline{\lambda}) = \underline{h}_0$
and hence   that (\ref{eqhconc}) is satisfied and that $\underline{\lambda}$ satisfies
equation (\ref{eqfpeqn}). The theorem is proved. \qed  

\begin{Cy}\label{cynew} Let $\underline{\lambda} \in \{0\} \times \mathbb{R}_+^{M-2} \times
\{0\}$ be a solution to equation (\ref{eqfpeqn}) where $\min_j p_j > 0$. Let
$\lambda^* = \max_{j \in \{2, \ldots, M-1\}}\lambda_j$. Then $\lambda^* <
+\infty$ and   for all $a \in [a_0,1)$, where $a_0 = \frac{\lambda^*}{1 +
\lambda^*}$, $\underline{q}(a) = \left(\frac{1}{a} - 1\right)\underline{\lambda}$ satisfies

\begin{equation}\label{eqqa} \underline{q}(a) = \left(\frac{1}{a} - 1 \right ){\cal F} \left(\underline{h}
\left(\underline{p}, \left(\frac{a}{1-a}\right)\underline{q}(a)\right)\right).\end{equation}
\end{Cy}

\paragraph{Proof} This follows directly from the preceding theorem. \qed

\paragraph{Proof of theorem~\ref{thnb}} Let $\mathbb{P}$ denote the probability measure with respect to both the random walk $X$ and the independent random time $\tau$. There is a solution to the problem if
and only if there is a $\underline{q} \in \{0\} \times [0,1]^{M-2} \times \{0\}$ such that
  the matrix $P(\underline{q})$
constructed according to  definition~\ref{defdtgd}, which is the one step
transition matrix for a discrete time Markov chain $X$,   satisfies 

\begin{eqnarray}\label{eqnbprob} \lefteqn{p_j =
\mathbb{P} (X_\tau = i_j|X_0 = e_0)}\\&& =  \begin{array}{ll}
(1-a)^r  \sum_{k=0}^\infty \left(\begin{array}{c} r+k-1 \\ r-1
\end{array}\right) \left(\frac{i_l - e_0}{i_l - i_{l-1}}((aP(\underline{q}))^k)_{l-1,j}  +
\frac{e_0 - i_{l-1}}{i_l - i_{l-1}}((aP(\underline{q}))^k)_{l,j} \right) & i_{l-1} < e_0 \leq i_l.  
                                              \end{array} \nonumber 
\end{eqnarray}

\noindent Define $G(\underline{q})$ by 

\begin{equation}\label{eqgedef} G(\underline{q}) =  \sum_{k=0}^\infty \left(\begin{array}{c} r
+ k - 1 \\ r -1\end{array}\right ) (aP(\underline{q}))^k,\end{equation}

\noindent then $G(\underline{q})$ satisfies
 
\begin{equation}\label{eqgee} G(\underline{q})(I-aP(\underline{q}))^r = (I - aP(\underline{q}))^r G(\underline{q}) = I.
\end{equation}

\noindent Using equation (\ref{eqgedef}),  $\underline{q} \in \{0\} \times (0,1]^{M-2} \times \{0\}$ is a solution to equation (\ref{eqnbprob}) if and
only if  

\begin{equation}\label{eqem} p_j =    
(1-a)^r\left(\frac{e_0 - i_{l-1}}{i_l - i_{l-1}}G_{l,j} (\underline{q}) + \frac{i_l -
e_0}{i_l - i_{l-1}}G_{l-1,j}(\underline{q})\right) \qquad  j = 1, \ldots, M \qquad  i_{l-1} < e_0 \leq 
i_l
\end{equation}
 
\noindent Recall the function $\underline{h}$ defined in equation (\ref{eqaitchemque})
satisfies equation (\ref{eqha}). It follows that 

\begin{eqnarray*} \frac{1}{1-a} \sum_{k=1}^M h_k \left (\underline{p},
\left(\frac{a}{1-a}\right) \underline{q} \right )(I - aP(\underline{q}))_{ k, j} &=&
\frac{1}{(1-a)^r} \sum_{k=1}^M p_k ((I - aP(\underline{q}))^r)_{ k, j}\\
 &=& \frac{1}{(1-a)^r}\sum_{k=1}^r p_k G^{-1}_{k, j}(\underline{q}).
\end{eqnarray*}

\noindent From equation (\ref{eqem}), using $h_k$ to denote $h_k\left(\underline{p}, \left(\frac{a}{1-a}\right)\underline{q} \right)$, it follows that $\underline{q} \in \{0\} \times (0,1]^{M-2} \times \{0\}$ is a
solution if and only if (using the expression in equation (\ref{eqnbprob}) for $\underline{p}$) 

\begin{eqnarray} \nonumber \sum_{k=1}^M h_k (I - aP(\underline{q}))_{ k, j}& = & 
\frac{1}{(1-a)^{r-1}} \sum_{k=1}^M p_k G^{-1}_{k,j}(\underline{q}) \\ &= & \nonumber (1-a)
\frac{e_0 - i_{l-1}}{i_l - i_{l-1}}\sum_{k=1}^M  G_{l,k}(\underline{q})G^{-1}_{k,j}(\underline{q}) + \frac{i_l - e_0}{i_l - i_{l-1}} \sum_{k=1}^M G_{l-1,k}(\underline{q}) G^{-1}_{k,j}(\underline{q}) \\ &=&  
\left\{\begin{array}{cc} \frac{i_l - e_0}{i_l - i_{l-1}}(1-a) & j = l-1 \\
\frac{e_0 - i_{l-1}}{i_l - i_{l-1}}(1-a) & j = l \\ 0 & j \in \{1,\ldots, M\}
\backslash \{l-1, l\} \end{array}\right.\label{eqhsoln}  \end{eqnarray}

\noindent if $i_{l-1} < e_0 \leq i_l$. From this, it follows that for any solution $\underline{q} \in \{0\} \times (0,1]^{M-2} \times \{0\}$,

\[ h_k =  
 (1-a) \left(\frac{i_l - e_0}{i_l - i_{l-1}} \sum_{j=0}^\infty a^j
(P(\underline{q})^j)_{i_{l-1}, i_k} + \frac{e_0 - i_{l-1}}{i_l - i_{l-1}} \sum_{j=0}^\infty
a^j (P(\underline{q})^j)_{i_l,i_k}\right) \qquad i_{l-1} <  e_0  \leq i_l \]

\noindent for $k = 1,\ldots, M$ and therefore, for $a \in (0,1)$,
for any solution $\underline{q} \in \{0\} \times (0,1]^{M-2} \times \{0\}$, $h_j$ is positive for each $j$.  From lemma~\ref{lmmaitchprex0}, it
follows that $\underline{h}$ is a  probability function on $S = \{i_1, \ldots, i_M\}$  such
that $\sum_{j=1}^M i_j h_j = e_0(\underline{p})$.\vspace{5mm}

\noindent It follows from lemma~\ref{lmmone} that $\underline{q}$ satisfying equation
(\ref{eqhsoln})  satisfies

\begin{equation}\label{eqqsol} q_j = \left(\frac{1}{a} - 1 \right) {\cal F}\left (\underline{h}
\left (\underline{p},\left(\frac{a}{1-a}\right) \underline{q} \right ), j \right ) \qquad j = 1, \ldots,
M.\end{equation} 

\noindent It follows from corollary~\ref{cynew} that there exists an  $a_0 < 1$
such that for all $a \in [a_0,1)$ there is a $\underline{q} \in \{0\} \times (0,1]^{M-2}
\times \{0\}$ satisfying equation (\ref{eqqsol}). The proof of
theorem~\ref{thnb} is complete. \qed 
 
\section{Proofs of theorems~\ref{thgdl1} and~\ref{thgdl2}}\label{appthgdl}

This section presents the proofs of theorems~\ref{thgdl1} and~\ref{thgdl2}. Theorem~\ref{thgdl1} states that for any prescribed law $\underline{p}$ over $S = \{i_1, \ldots, i_M\}$ and $T \sim \Gamma \left( r, \frac{t}{r} \right)$ a random time with a gamma distribution, there is a CMRW (definition~\ref{defctgd}) $X$, such that for $T$ independent of $X$,  $\mathbb{P}(X_T = i_j | X_0 = e_0(\underline{p})) = p_j$ for $j = 1, \ldots, M$. Theorem~\ref{thgdl2} takes the limit as $r \rightarrow +\infty$, with $\mathbb{E}[T] = t$ fixed. It follows that $T \rightarrow t$ in law, from which it follows that for any $t > 0$, there is a CMRW $X$ such that $\mathbb{P}(X_t = i_j |X_0 = e_0(\underline{p})) = p_j$ for $j = 1, \ldots, M$, the subject of theorem~\ref{thgdl2}.   

Theorem~\ref{thgdl1} is proved using theorem~\ref{thnb}. Let $\tau \sim NB(r,a)$ and $T^{(\delta)} = \tau \delta$. For fixed $r$, the parameters $a$ and $\delta$ are chosen such that $\mathbb{E}[T^{(\delta)}] = t$. As $\delta \rightarrow 0$ and $a \uparrow 1$, $T^{(\delta)} \stackrel{(d)}{\longrightarrow} \Gamma \left( r, \frac{t}{r} \right )$. 

Subsection~\ref{appthgdlsub1} deals with the proof of theorem~\ref{thgdl1}, while subsection~\ref{appthgdlsub2} deals with the proof of theorem~\ref{thgdl2}.  

\subsection{Proof of theorem~\ref{thgdl1}}\label{appthgdlsub1} 
By theorem~\ref{thnb}, there exists an
$a_0$ such that for all $a \in [a_0,1)$, if $\tau \sim NB(r,a)$ (negative
binomial with parameters $r$ and $a$), then there exists a DMRW $X$ (depending on the parameter $a$), independent of $\tau$, such that 

\[ \mathbb{P} \left ( X_\tau = i_j | X_0 = e_0(\underline{p}) \right ) = p_j.\]

\noindent The one step transition probability matrix of $X$ is $P(\underline{q})$ (following
definition~\ref{defdtgd}), where $\underline{q} = \left( \frac{1}{a}  - 1 \right )\underline{\lambda}$,
where $\underline{\lambda}$ is a fixed point of equation (\ref{eqfpeqn}). 
For fixed $r \geq 1$, let $T^{(\delta)} = \tau \delta$ and choose $\delta$ such
that $\mathbb{E} [T^{(\delta)}] = t$. Then, since $\mathbb{E} [T^{(\delta)}] = \frac{ra\delta}{(1-a)} = t$, it follows that $a = \frac{1}{1 + \left(\frac{r \delta}{t}\right)}$ giving $\underline{q} = \frac{r\delta}{t} \underline{\lambda}$.

\noindent Choose $\delta$ sufficiently small such that $a > a_0$; that is 
\[ \delta < \delta_0 = \frac{t}{r}\left(\frac{1}{a_0} - 1\right).\]

\begin{Lmm} \label{lmmnbgam} Let ${\cal L}(T^{(\delta)})$ denote the law of $T^{(\delta)}$. Then, as $\delta \rightarrow 0$,
\[ {\cal L}(T^{(\delta)}) \stackrel{\delta \rightarrow 0}{\longrightarrow}
\Gamma \left (r, \frac{t}{r} \right ).\]
\end{Lmm}

\noindent That is, let $T$ denote a random variable with the limiting
distribution, then $T$ has a Gamma distribution. This may be expressed as: $T =
X_1 + \ldots, X_r$, where $X_j \sim Exp(\frac{r}{t})$ and $X_1, \ldots, X_r$ are
independent.

\paragraph{Proof of lemma~\ref{lmmnbgam} } Firstly,  
\[ \mathbb{P}(T^{(\delta)} = n\delta) = \left(\begin{array}{c} n+r-1 \\ r-1
\end{array}\right) \left(\frac{r\delta}{t}\right)^r \frac{1}{(1 +
\frac{r\delta}{t})^{n+r}} \qquad n = 0, 1, 2, \ldots \]

\noindent and, for $T \sim \Gamma \left(r, \frac{t}{r} \right)$, 

\begin{eqnarray*} \mathbb{P}(n\delta \leq T < (n+1)\delta) & = & \frac{1}{
(r-1)!}\left(\frac{r}{t}\right)^r \int_{n\delta}^{(n+1)\delta} x^{r-1}e^{-rx/t}
dx \\
&=& \frac{1}{(r-1)!}\int_{rn\delta/t}^{r(n+1)\delta/t} y^{r-1}e^{-y} dy
\end{eqnarray*}

\noindent It follows that 

\[  \frac{r\delta }{t (r-1)!} \min_{y \in \left [\frac{rn\delta}{t},
\frac{r(n+1)\delta}{t}\right ]} y^{r-1}e^{-y} \leq \mathbb{P}(n\delta \leq T <
(n+1)\delta ) \leq  \frac{r\delta }{t (r-1)!} \max_{y \in \left
[\frac{rn\delta}{t}, \frac{r(n+1)\delta}{t}\right ]} y^{r-1}e^{-y} \]

\noindent from which it follows that for any $s > 0$, with $n = [s/\delta]$,
using Stirling's formula, 

\begin{eqnarray*}\lefteqn{ \lim_{\delta \rightarrow 0}   \left | \frac{\mathbb{P} (T^{(\delta)} = 
n\delta)}{\mathbb{P} (n \delta \leq T < \left(n + 1\right)\delta)}  - 1 \right |  = \lim_{\delta
\rightarrow 0}   \left | \frac{ \frac{(n+r-1)!}{(r-1)!n!}\left(\frac{r\delta}{t}\right)^r
(1 + \frac{r\delta}{t})^{-(n+r)}}{\frac{1}{(r-1)!}\left(\frac{r}{t}\right)^{r}
s^{r-1}\delta e^{-rs/t}} - 1 \right |} \\&& = \lim_{\delta \rightarrow 0}   \left | \frac{e^{-(r-1)}(1 +
\frac{r-1}{n})^{n+r-(1/2)} n^{r-1} \delta^{r-1}(1+
\frac{r\delta}{t})^{-(n+r)}}{s^{r-1}e^{-rs/t}}  - 1 \right | \\
&& = \lim_{\delta \rightarrow 0}   \left | \frac{e^{-(r-1)}(1 + \frac{r-1}{n})^{n+r-(1/2)}
 (1+ \frac{r\delta}{t})^{-(n+r)}}{ e^{-rs/t}} - 1 \right |   = 0. 
 \end{eqnarray*}

\noindent From this, it follows that

\[ \lim_{\delta \rightarrow 0} \sup_{s \geq 0} \left | \mathbb{P} (T^{(\delta)} \leq s) - \mathbb{P}(T \leq s) \right | = 0\] and the result follows.  \qed \vspace{5mm}

\noindent Let $Y^{(\delta)}_s = X_{[s / \delta]}$. Convergence in distribution in the sense of  finite dimensional marginals is outlined. That is, 
for any $0 < s_1 < \ldots < s_n \leq t$ and any $(i_{k_0}, \ldots, i_{k_n})
\subset S^n$, convergence of 

\[ \mathbb{P} \left (Y^{(\delta)}_{s_1} = i_{k_1}, \ldots, Y^{(\delta)}_{s_n} = i_{k_n} |
Y^{(\delta)}_{0} = i_{k_0} \right ) = \prod_{j=1}^n \mathbb{P} \left (Y^{(\delta)}_{s_j}
= i_{k_j} | Y^{(\delta)}_{s_{j-1}} = i_{k_{j-1}} \right ).\]

\noindent is outlined. Note that

\[\mathbb{P} \left (Y_{s_j}^{(\delta)} = i_{k_j}\left | Y_{s_{j-1}}^{(\delta)} =
i_{k_{j-1}} \right. \right ) = \mathbb{P} \left
(Y^{(\delta)}_{\delta\left[\frac{s_j}{\delta}\right]} = i_{k_j} \left |
Y^{(\delta)}_{\delta \left [ \frac{s_{j-1}}{\delta}\right ]} =
i_{k_{j-1}}\right. \right ).\]

\noindent For any $0 < s_1 < s_2 < +\infty$,

\[\mathbb{P} \left( \left. Y^{(\delta)}_{s_1 + r} = i_j \quad \forall \quad 0 < r < s_2
\quad \right | \quad Y^{(\delta)}_{s_1} = i_j \right) = \left ( 1 -
q_j \right)^{f(s_2)} =  \left( 1 -
\frac{r\delta}{t}\lambda_j \right)^{f(s_2)}
\]

\noindent where $f(s_2)$ is either $[s_2/\delta]$ or $[s_2/\delta] + 1$. It
follows that 

\begin{equation}\label{eqtaulimcdf} \mathbb{P} \left(Y^{(\delta)}_{s_1 + r} = i_j \quad
\forall \quad 0 < r < s_2 \quad \left | \quad  Y^{(\delta)}_{s_1} = i_j \right.
\right) \stackrel{\delta \rightarrow 0}{\longrightarrow} \exp\left\{
-\left(\frac{r}{t}\lambda_j \right)s_2\right\}.\end{equation}

\noindent Let $\tau^{(\delta)}(i_j)$ denote a random variable with cumulative
distribution function

\begin{equation}\label{eqtaucdf}  \mathbb{P} (\tau^{(\delta)}(i_j) \leq s) =  \left\{ \begin{array}{ll} 1 - \left( 1 -
\frac{r\delta}{t} \lambda_j \right)^{[s/\delta]} & \delta > 0 \\
 1 -  \exp \left \{ - \frac{rs}{t} \lambda_j \right \} & \delta = 0 \end{array}\right.\end{equation}

\noindent and let \[ \tau^{(\delta)}_0(i_{j_0}), \tau^{(\delta)}_1(i_{j_1}),
\tau^{(\delta)}_2(i_{j_2}) , \ldots\] 

\noindent denote independent random variables such that $\tau^{(\delta)}_n( i_j)
\stackrel{(d)}{=}\tau^{(\delta)}(i_j)$ for all $n \geq 1$ and all $i_j \in S$. Then

\begin{eqnarray}\lefteqn{ \mathbb{P} \left
(Y^{(\delta)}_{\delta\left[\frac{s}{\delta}\right]} = i_{j} \left |
Y^{(\delta)}_{0} = i_{k}\right. \right ) = \mathbb{P}(\tau_0^{(\delta)}(i_k) > s; j =
k)} \label{eqfdmcon} \\&& + \sum_{n=1}^\infty \sum_{a_1,\ldots, a_{n-1}} \mathbb{P}
\left (\tau_0^{(\delta)} (i_k) + \sum_{m=1}^{n-1} \tau^{(\delta)}_m(i_m)   < s <
\tau_0^{(\delta)} (i_k) + \sum_{n=1}^{n-1} \tau^{(\delta)}_m(i_m)   +
\tau^{(\delta)}_n(i_j) \right )\nonumber \\&& \hspace{110mm} \times \mathbb{P}( 
a_1,\ldots,a_{n-1},i_j)\nonumber 
\end{eqnarray}

\noindent where the notation $\mathbb{P}(a_1,\ldots,a_{n-1},i_j)$ denotes the probability
that at the first $n$ jump times, the sequence of sites the process visited was
$a_1,\ldots, a_{n-1}, i_j$. It follows from equation (\ref{eqtaucdf}) (details omitted) that  

\begin{eqnarray*}\lefteqn{\lim_{\delta \rightarrow 0} \sup_{s > 0} \max_{(j,k) \in \{1, \ldots, M\}^2} \left | \mathbb{P} \left
(\tau_0^{(\delta)} (i_k) + \sum_{m=1}^{n-1} \tau^{(\delta)}_m(a_m)   < s <
\tau_0^{(\delta)} (i_k) + \sum_{n=1}^{n-1} \tau^{(\delta)}_m(a_m)   +
\tau^{(\delta)}_n(i_j) \right ) \right.} \\&& \left. - \mathbb{P} \left (\tau_0^{(0)}
(i_k) + \sum_{m=1}^{n-1} \tau^{(0)}_m(a_m)   < s < \tau^{(0)}_0 (i_k) +
\sum_{n=1}^{n-1} \tau^{(0)}_m(a_m)   + \tau_n^{(0)}(i_j) \right ) \right | =
0\end{eqnarray*}

\noindent for each $n \geq 1$ and that 

\[ \lim_{\delta \rightarrow 0} \sup_{s > 0} \max_{(j,k) \in \{1, \ldots, M\}}\left | \mathbb{P}(\tau_0^{(\delta)}(i_k) > s; j= k) -
\mathbb{P}(\tau_0^{(0)}(i_k) > s ; j = k)\right | = 0.\]

\noindent From this, it follows (details omitted) that the processes $Y^{(\delta)}$ converge, in the sense of  convergence of finite dimensional marginals, to a process $Y$. In particular, $Y$ satisfies  

\begin{equation}\label{eqsupconv} \lim_{\delta \rightarrow 0} \sup_{s > 0} \max_{(j,k) \in \{1, \ldots, M\}^2} \left |  \mathbb{P} \left
(Y^{(\delta)}_{\delta\left[\frac{s}{\delta}\right]} = i_{j} \left |
Y^{(\delta)}_{0} = i_{k}\right. \right ) - \mathbb{P} \left (Y_s = i_{j}  \left |
Y_{0} = i_{k} \right. \right )   \right |  = 0.
\end{equation}

\noindent Note that, for all $\delta \in (0, \delta_0)$,

\[ \mathbb{P}(Y^{(\delta)}_{T^{(\delta)}} = i_j | Y^{(\delta)}_0 = e_0) = p_j.\]

\noindent  It is required to show that this also holds in the limit; namely,
$\mathbb{P}(Y_T = i_j | Y_0 = e_0) = p_j$. Firstly,

\begin{eqnarray*}\lefteqn{ |\mathbb{P}(Y_{T^{(\delta)}}^{(\delta)} = i_j | Y_0^{(\delta)}
= e_0) - \mathbb{P}(Y_T = i_j | Y_0 = e_0)|  \leq |\mathbb{P}(Y^{(\delta)}_{T^{(\delta)}} =
i_j | Y_0^{(\delta)} = e_0) - \mathbb{P}(Y_{T^{(\delta)}} = i_j|Y_0 = e_0)|}\\&& \hspace{65mm}  +
|\mathbb{P}(Y_{T^{(\delta)}} = i_j | Y_0 = e_0) - \mathbb{P}(Y_T = i_j|Y_0 = e_0)|  = I + II.
\end{eqnarray*}

\noindent For $I$,

\begin{eqnarray*}
\lefteqn{ |\mathbb{P}(Y^{(\delta)}_{T^{(\delta)}} = i_j | Y_0^{(\delta)} = e_0) -
\mathbb{P}(Y_{T^{(\delta)}} = i_j|Y_0 = e_0)|}\\&&
\leq \sum_{n=0}^\infty \left | \mathbb{P}(Y_{n\delta}^{(\delta)} = i_j | Y_0^{(\delta)} =
e_0) - \mathbb{P}(Y_{n\delta} = i_j | Y_0 = e_0)\right |\mathbb{P}(T^{(\delta)} = n\delta)\\&&
\leq \sum_{n=0}^\infty \left | \mathbb{P}(Y_{n\delta}^{(\delta)} = i_j | Y_0^{(\delta)} =
e_0) - \mathbb{P}(Y_{n\delta} = i_j | Y_0 = e_0)\right |\frac{\mathbb{P}(T^{(\delta)} = n\delta)}{\mathbb{P}(n\delta \leq T < (n+1)\delta)}\mathbb{P}(n\delta \leq T < (n+1)\delta)\\&&
= \int_0^\infty g_n(x) f_T(x) dx
\end{eqnarray*}
\noindent where $f_T$ is the density function of $T$ and 
\[ g_n(x) =   \left | \mathbb{P}(Y_{n\delta}^{(\delta)} = i_j | Y_0^{(\delta)} =
e_0) - \mathbb{P}(Y_{n\delta} = i_j | Y_0 = e_0)\right |\frac{\mathbb{P}(T^{(\delta)} = n\delta)}{\mathbb{P}(n\delta \leq T < (n+1)\delta)} \qquad n\delta \leq x < (n+1)\delta.     \]

\noindent It now follows directly from equation (\ref{eqsupconv}) and  dominated convergence that $I \stackrel{\delta \rightarrow 0}{\longrightarrow} 0$. 
 \vspace{5mm}

\noindent For $II$, let $f_T$ denote the probability density function
of $T$, then 

\begin{eqnarray}\nonumber 
\lefteqn{ |\mathbb{P}(Y_{T^{(\delta)}} = i_j | Y_0  = e_0) - \mathbb{P}(Y_{T} = i_j|Y_0 =
e_0)|}\\&& \nonumber 
\leq \sum_{n=0}^\infty \left | \mathbb{P}(Y_{n\delta}  = i_j | Y_0^{(\delta)} =
e_0)\mathbb{P}(T^{(\delta)} = n\delta) - \int_{n\delta}^{(n+1)\delta} \mathbb{P}(Y_t = i_j | Y_0 =
e_0)f_T(t) dt \right |\\&& \nonumber 
\leq \sum_{n=0}^\infty   \mathbb{P}(Y_{n\delta}  = i_j | Y_0^{(\delta)} = e_0)\left |
(\mathbb{P}(T^{(\delta)} = n\delta) - \mathbb{P}(n\delta \leq T < (n+1)\delta)) \right|\\&&
\hspace{5mm} +  \sum_{n=0}^\infty  \int_{n\delta}^{(n+1)\delta} \left |\mathbb{P} (Y_{n\delta} = i_j | Y_0 = e_0) -  \mathbb{P} (Y_t = i_j |
Y_0 = e_0) \right | f_T(t) dt \nonumber\\ && = II(i) + II(ii).\label{eqtconv}
\end{eqnarray}

\noindent For $II(i)$, 

\begin{eqnarray*} \lefteqn{\sum_{n=0}^\infty   \mathbb{P}(Y_{n\delta}  = i_j | Y_0^{(\delta)} = e_0)\left |
(\mathbb{P}(T^{(\delta)} = n\delta) - \mathbb{P}(n\delta \leq T < (n+1)\delta)) \right|}\\&& \leq 
\sum_{n=0}^\infty \mathbb{P} \left (n\delta \leq T < (n + 1) \delta \right ) \left | 1 - \frac{\mathbb{P}(T^{(\delta)} = n\delta)  }{ \mathbb{P} (n \delta \leq T < (n+1)\delta)} \right |  \stackrel{\delta \rightarrow 0}{\longrightarrow} 0
\end{eqnarray*}

\noindent For $II(ii)$,  

\[ \max_n \max_{t \in [n\delta, (n+1)\delta]}|\mathbb{P}(Y_{n\delta} = i_j |Y_0 = e_0) - \mathbb{P}(Y_t =
i_j|Y_0 = e_0)| \stackrel{\delta \rightarrow 0}{\longrightarrow 0}\]

\noindent from which it follows that $II(ii) \stackrel{\delta \rightarrow 0}{\longrightarrow} 0$.  \vspace{5mm} 

 \noindent The proof of theorem~\ref{thgdl1} is complete. \qed 

\subsection{Proof of theorem~\ref{thgdl2}}\label{appthgdlsub2}  By theorem~\ref{thgdl1}, for any integer $r \geq 1$, any given $t > 0$ and any given fixed probability measure $\underline{p}$ over a finite state space $S = \{i_1, \ldots, i_M\}$, there exists a CMRW (definition~\ref{defctgd}) $X^{(r)}$ such that $\mathbb{P} \left (X_T^{(r)} = i_j | X_0^{(r)} = e_0(\underline{p}) \right ) = p_j$, where $T$ is independent of $X^{(r)}$ and satisfies  $T \sim \Gamma(r, \frac{t}{r})$. Let $\underline{h}^{(r)}$ satisfy $\underline{h}^{(r)}(\underline{p}, \underline{\lambda}) = \underline{p} N^{r-1}(\underline{\lambda})$. This is equation (\ref{eqaitchemque}), where the parameter $r$ is explicitly stated, because the task here is to let $r \rightarrow +\infty$.  Let $\underline{\lambda}^{(r)}$ denote the fixed point, satisfying equation (\ref{eqfpeqn}) from theorem~\ref{thfpt}. From the proof of theorem~\ref{thgdl1}, the Markov chain $X^{(r)}$ remains on a $i_j$ for an exponential holding time $\tau^{(r)} (i_j)$ which satisfies $\mathbb{P} \left (\tau^{(r)} (i_j) \geq s \right ) = \exp\left\{ - \frac{rs}{t} \lambda_j^{(r)} \right\}$. This is the content of equation (\ref{eqtaucdf}).  Now set $\rho_j^{(r)} = \frac{r}{t} \lambda_j^{(r)}$ and let $\rho^{(r)*} = \max_{j \in \{2, \ldots, M-1\}} \rho_j^{(r)}$. Suppose that $\sup_r \rho^{(r)*} = +\infty$. It follows that there is a $j \in \{2, \ldots, M-1\}$ and a sequence $(r_n)_{n \geq 1}$ such that $\rho_j^{(r_n)} \rightarrow +\infty$. Then, for $i_j \geq e_0$ and $k \in \{j+1, \ldots, M\}$, or $i_j \leq e_0$ and $k \in \{1,\ldots,j-1\}$, 

\[ p_{j+1}  \leq 1 - \mathbb{P} (\tau^{(r_n)}(i_j) > T) = 1 - \int_0^\infty  \mathbb{P} (\tau^{(r_n)}(i_j) \geq x) f_T(x) dx\stackrel{n \rightarrow +\infty}{\longrightarrow} 0 \]

\noindent where $f_T$ denotes the density function of a $\Gamma \left( r, \frac{t}{r} \right)$ random variable. This contradicts the fact that $\min_{j \in \{1, \ldots, M\}} p_j > 0$.   Therefore,   $\sup_r \rho^{(r)*} < +\infty$ and therefore there exists a limit point $\underline{\rho}$ of $(\underline{\rho}^{(r)})_{r \geq 0}$. 

Let $T_r \sim \Gamma \left(r, \frac{t}{r} \right)$. Then $\mathbb{E} [T_r] = t$ and $\mathbb{V}(T_r) = \frac{t^2}{r} \stackrel{r \rightarrow +\infty}{\longrightarrow} 0$, so that for any $\epsilon > 0$, it follows from Chebychev's inequality that 
\[ \mathbb{P}\left ( \left |  T_r - t\right | > \epsilon\right ) \stackrel{r \rightarrow +\infty}{\longrightarrow} 0.\] 

Let $(\underline{\rho}^{(r_n)})_{n \geq 1}$ denote a subsequence such that $\lim_{n \rightarrow +\infty} \max_{j \in \{1, \ldots, M\}} |\rho_j - \rho_j^{(r_n)}| = 0$. Let $X$ denote the CMRW process on $S = \{i_1, \ldots, i_M\}$ with holding intensities $\underline{\rho}$. Then for each $j$ and each $r$, $\mathbb{P}(X_{T_r}^{(r)} = i_j) = p_j$. Furthermore,

\[ \lim_{n \rightarrow +\infty} \max_{j \in \{1, \ldots, M\}} |\mathbb{P}(X_{T_{r_n}}^{(r_n)} = i_j) - \mathbb{P}(X_t^{(r_n)} = i_j)| = 0\]

\noindent and

\[ \lim_{n \rightarrow +\infty} \max_{j \in \{1, \ldots, M\}} |\mathbb{P}(X_t^{(r_n)} = i_j) - \mathbb{P}(X_t = i_j)| = 0\]

\noindent so that 

\[ \mathbb{P}(X_t = i_j) = p_j \qquad j = 1, \ldots, M.\]

\noindent The theorem is proved.  \qed 

\section{Kre\u{i}n's Strings and Markov processes}\label{seckrstr} 

The discussion so far has been about finite state spaces and constructing Markov processes, both discrete and continuous time, that have a given marginal at an independent random time with a certain specified distribution (geometric, negative binomial, gamma), then taking a limit to obtain a specified deterministic time, with the requirement that the process is a martingale, which jumps only to nearest neighbours.

The remainder of the discussion is the extension to prove existence of a Markov martingale of diffusion type with specified distribution at a specified time $t > 0$, for an arbitrary specified distribution $\mu$, such that $\int_{-\infty}^\infty |x| \mu(dx) < +\infty$, over $\mathbb{R}$. This is the content of theorem~\ref{thgdcont}. 

The process in question, in the generality required, is discussed by Kotani and Watanabe in~\cite{KoWa}, where it is called a `generalised diffusion' and by Knight in~\cite{Kn}, where it is called a `gap diffusion'. The term `skip free diffusion' is also used in the literature, emphasising that the movement is between nearest neighbours. 

Section~\ref{seckrstr} gives the background on Kre\u{i}n's strings, their relation to Markov processes and the results on spectral theory that are needed for some technical lemmas. Section~\ref{appthgdcont} then considers a limiting sequence of atomised measures and shows that a subsequence of the corresponding processes converge to a well defined gap diffusion with the prescribed terminal distribution. 

\subsection{Kre\u{i}n's strings} Kre\u{i}n defined a string  on an interval $[0, l]$, for $l \leq +\infty$, where $l$ is the {\em length} of the string and, for any Borel set $A \in {\cal B}([0,l])$, $m(A)$ corresponds to the mass of the string on set $A$. The string is `tied down' at the point $0$. Kre\u{i}n was interested in the dynamics of the string. He introduced the second order operator $\frac{d^2}{dm dx}$ (defined formally later) with this in view. He developed the theory of such operators in~\cite{Kr}; their
spectral properties were studied in~\cite{KK}. They are discussed at length in
the book~\cite{DM}. The reader is referred to these works for proofs of results
about strings that are stated here without proof.

This operator may be seen as the infinitesimal generator of a Markov process, which Kre\u in did not have in view. To give a full description of the Markov process, Kre\u in's operator was altered to accommodate Dirichlet boundary conditions at the left end point,  as found in Kotani and  Watanabe~\cite{KoWa}.  

\paragraph{Notation} Throughout, the following notations will be used:
\begin{enumerate}
 \item For an interval, a square bracket denotes that the end point is included;
a curved bracket denotes that it is not included. For $x < y$, $(x,y)$ includes
neither $x$ nor $y$, $[x,y]$ includes both $x$ and $y$, $[x,y)$ includes $x$ but
not $y$ and $(x,y]$ includes $y$ but not $x$.

\item For $x \in \mathbb{R}$, $[x]$ denotes the integer part of $x$. That is, 
 \[ [x] = \max \{z \in \mathbb{Z} | z \leq x\}.\]
 \item $f^\prime_-(x)$ denotes the left derivative; $\lim_{h \rightarrow
0}\frac{f(x) - f(x-h)}{h}$ and $f^\prime_+(x)$ denotes the right derivative
$f^\prime_+(x) = \lim_{h \rightarrow 0} \frac{f(x+h) - f(x)}{h}$.
 \item  $\int_a^b$ signifies integration over the half closed interval $a < x
\leq b$. Similarly, $\int_{a-}^b$ indicates that the mass $x = a$ is included
while $\int_{a-}^{b-}$ includes the mass of $a$ but excludes the mass of $b$.
\end{enumerate}

\begin{Defn}[Kre\u{i}n String] \label{defks} A {\em string} measure $m^*$ is defined as a measure over ${\cal B}({\bf R})$,  that satisfies the following properties: there is a   non-decreasing càdlàg (right continuous, left limits) function,  ${\cal M}$, defined on $(-\infty, +\infty)$  such that ${\cal M}(0) = 0$ and for $x > -\infty$, ${\cal M}(x) > -\infty$ and for $x < +\infty$, ${\cal M}(x) < +\infty$.  Let  $l_0 = \sup\{x | {\cal M} (x) = \inf_z {\cal M} (z)\}$ and $l_1 =   \inf\{x|{\cal M}(x) = \sup_z {\cal M}(z)\}$. If $l_1 < +\infty$, there is a point $L_1 \in [l_1,+\infty)$ and if $l_0 > -\infty$ there is a point $L_0 \in (-\infty, l_0]$ ($L_1 = +\infty$ if $l_1 = +\infty$ and $L_0 = -\infty$ if $l_0 = -\infty$)  
such that 
 \begin{equation} \label{eqmfromM2} \left\{ \begin{array}{l} m^*((a,b]) =  {\cal M}(b) - {\cal M} (a) \qquad L_0 \leq a \leq b < L_1 \\  m^*(\{L_1\}) = +\infty \quad \mbox{and} \quad m^*((L_1,+\infty)) = 0 \quad \mbox{if}\quad L_1 < +\infty\\
 m^*(\{L_0\}) = +\infty \quad \mbox{and} \quad m^*((-\infty, L_0)) = 0 \quad \mbox{if} \quad L_0 > -\infty \end{array}\right.  \end{equation}
  \end{Defn}\vspace{5mm}
  
\noindent  $m^*$ is used for construction of the Markov process according to Kotani and Watanabe~\cite{KoWa} with cemetaries $L_0$ and $L_1$.

The gap diffusion process of definition~\ref{defmpsting} has infinitesimal generator ${\cal G}^* = \frac{d^2}{dm^* dx}$. As S. Kotani points out and proves in  appendix 1 of~\cite{KoWa}, pp 250 - 254, the operator ${\cal G}^*$ defines a self adoint, positive semidefinite operator on the Hilbert space $L^2([L_0, L_1], m^*)$; the difference between the string of Kac and Kre\u in and $m^*$ discussed by Kotani is the difference between Neuman boundary conditions at $l_0$ and Dirichlet boundary conditions at $L_0$.

\paragraph{Notation} For $A \in {\cal B}(\mathbb{R})$, 
 \[ \|f\|_{L^2(A,m^*)}^2 = \int_A f^2(x)m^*(dx)\qquad \mbox{and} \qquad \langle
f,g \rangle_{L^2(A,m^*)} = \int_A f(x)g(x) m^*(dx) \] 
with the convention that if $h(x) = 0$ and $m^*(\{x\}) = +\infty$, then $h(x)m^*(\{x\}) = 0$. 
 
\subsection{The differential operator}  The differential operator ${\cal G}^*$ associated with a string $m^*$ is  defined by

\[{\cal G}^* f = \frac{d}{dm^*} f^\prime_+,\]

\noindent where $f^\prime_+$ and $f^\prime_-$ denote right derivatives and left
derivatives respectively. The formal definition is taken `under the integral'
and is given in equation (\ref{eqdiffop}).

When $L_1 - L_0 < +\infty$, the operator ${\cal G}^* = \frac{d^2}{dm^* dx}$, restricted to ${\cal D}({\cal G}^*)$, the domain of the
operator (defined below), is a self adjoint, negative definite, densely defined operator on the
Hilbert space $L^2([L_0,L_1],m^*)$. The relevant results from~\cite{Kr}, \cite{KK} and \cite{DM} with modifications to the Dirichlet boundary conditions at $L_0$ given by Kotani in appendix 1 to~\cite{KoWa}, are summarised below.  

\begin{Defn} [Operator domains: $L_1 - L_0 < +\infty$, $m((L_0,L_1)) < +\infty$] \label{defopdom} Let $m^*$ be a string defined from a function ${\cal M}$ using  equation (\ref{eqmfromM2}), with $L_1 - L_0 < +\infty$ and $m^*((L_0,L_1)) < +\infty$, $L_0$ and $L_1$ are as in definition~\ref{defks}. 
The domain ${\cal D}_0({\cal G}^*)$ is defined as the space of functions $f$
defined on the whole real line $\mathbb{R}$ that satisfy the following
property: there is a function $g$ such that

\begin{equation}\label{eqdiffop} \left\{\begin{array}{ll} f(x) = f(L_0) + f^\prime_- (L_0) (x - L_0) & x \leq L_0 \\  f(x) = f(L_0)    +
f^\prime_-(L_0)(x - L_0) + \int_{L_0}^x \int_{L_0-}^y g(z)  m^*(dz) dy & x \geq L_0 \\
\int_{L_0-}^{L_1} g^2(z) m^*(dz) < +\infty & \end{array}\right.  
\end{equation}

\noindent Such a function $g$ is written $g = {\cal G}^* f$. \vspace{5mm} 

\noindent The domain ${\cal D}_-({\cal G}^*)$ is defined as
\begin{equation}\label{eqdmin} {\cal D}_-({\cal G}^*) = {\cal D}_0({\cal G}^*) \cap \{f | f (x) = 0 \quad \forall x \leq L_0\}.\end{equation}

\noindent The domain ${\cal D}_+({\cal G}^*)$ is defined as

\begin{equation}\label{eqdplus} {\cal D}_+({\cal G}^*) = {\cal D}_0({\cal G}^*) \cap \{f | \|f\|_{L^2([L_0,L_1],m^*)} + \|{\cal G}^*f\|_{L^2([L_0,L_1],m^*)} <
+\infty\} \cap \{f(x) = 0 \quad \forall x \geq L_1\}.\end{equation}

\noindent The domain ${\cal D} ({\cal G}^*)$ is defined as
\begin{equation}\label{eqdomd} {\cal D} ({\cal G}^*) = {\cal D}_-({\cal G}^*) \cap {\cal D}_+({\cal G}^*).\end{equation}
 \end{Defn}
 
 \noindent  This implies that for $L_1 - L_0 < +\infty$, $m^*((L_0,L_1)) < +\infty$ and $m^*(\{L_0\}) = m^*(\{L_1\}) = +\infty$, any $f \in {\cal D}({\cal G}^*)$ satisfies $f(x) = 0
\quad \forall x  \in (-\infty, L_0] \cup [L_1, +\infty)$ and  ${\cal G}^*f(L_0) = {\cal
G}^*f(L_1) = 0$; if ${\cal G}^*f(L_1) \neq 0$, then  $m^*(\{L_1\})({\cal G}^*f(L_1))^2 = +\infty$ and hence $\|{\cal G}^*f\|_{L^2([L_0,L_1],m^*)} = +\infty$, similarly for ${\cal G}^*f(L_0)$.    

\begin{Defn}[Operator domains: $L_1 - L_0 = +\infty$]\label{defopdom2}
 Let $m^*$ be a string defined from a function ${\cal M}$ using equation (\ref{eqmfromM2}) where $L_1 - L_0 = +\infty$.  The domain ${\cal D}_0({\cal G}^*)$ is defined as the space of functions $f$ defined on the whole real line $\mathbb{R}$ that satisfy the following property: if $L_0 > -\infty$, there is a function $g$ that satisfies equation (\ref{eqdiffop}). If $L_0 = -\infty$, then there is a function $g$ such that 
 
 \begin{equation}\label{eqdiffop2} \left\{\begin{array}{l} f(x) =  
  f(a  ) + (x- a) f^\prime_-(a ) + \int_{a }^x \int_{a-}^y g(z)  m^*(dz) dy     \\ \hspace{10mm} \forall -\infty < a < x < +\infty \\
\int_{x-}^{y} g^2(z) m^*(dz) < +\infty  \qquad \forall -\infty < x < y < +\infty.  \end{array}\right.  
\end{equation}

\noindent  Such a function $g$ is written $g = {\cal G}^*f$.  If $L_0 > -\infty$, then ${\cal D}_-({\cal G}^*)$ is defined by equation (\ref{eqdmin}). If $L_0 = -\infty$, then ${\cal D}_-({\cal G}^*)$ is defined as 
\[ {\cal D}_-({\cal G}^*) = {\cal D}_0({\cal G}^*) \cap \{f | \lim_{a \rightarrow +\infty} f^\prime_-(a) = 0, \quad \lim_{a \rightarrow +\infty} f(a) = 0\}.\]
If $L_1 - L_0 =  +\infty$, then 

\begin{equation}\label{eqdplus2} \left\{ \begin{array}{l} {\cal D}_+({\cal G}^*) = {\cal D}_0 ({\cal G}^*) \cap \{f | \|f\|_{L^2([z_1,z_2],m^*)} + \|{\cal G}^* f \|_{L^2([z_1,z_2],m^*)} < + \infty\} \\ \hspace{30mm} \cap \{f| \lim_{b \rightarrow +\infty} f(b) = 0,\; \lim_{b \rightarrow +\infty} f^\prime_+(b) = 0\} \\
\forall -\infty < z_1 < z_2 < +\infty, \quad z_1 \geq L_0, \quad z_2 \leq L_1. 
           \end{array} \right. \end{equation} 
 The domain ${\cal D}({\cal G}^*)$ is defined as ${\cal D}({\cal G}^*) = {\cal D}_-({\cal G}^*) \cap {\cal D}_+({\cal G}^*)$. 

\end{Defn}

\noindent The following theorem, stated without proof, describes the spectral
theory that will be used.

\begin{Th}   Let $m^*$ be a string defined by equation (\ref{eqmfromM2}) such that $L_1 - L_0 < +\infty$ and $m^*((L_0,L_1)) < +\infty$.  The set ${\cal D}({\cal G}^*)$ is a dense subset of
$L^2([L_0,L_1],m^*)$. Every $f \in {\cal D}({\cal G}^*)$ has left limits
and right derivatives, satisfies  $f_-^\prime(L_0) = f_+^\prime (L_1) = 0$, $f(x) = 0$ for $x \in (-\infty, L_0] \cup [L_1, +\infty)$ and
the operator ${\cal G}^* : {\cal D}({\cal G}^*) \rightarrow L^2([L_0,L_1],m^*)$
defined by ${\cal G}^*f = \frac{d}{dm^*} f^\prime_+$ is well defined, self adjoint
and non negative definite. Let $(\eta_k, -\lambda_k)_{k \geq 0}$ denote the
sequence of eigenfunctions / eigenvalues. Then $\inf_k \lambda_k > 0$. 
 \end{Th}

\paragraph{Proof} The result, including the strictly positive lower bound on the spectrum under the hypotheses that $L_1 - L_0 < +\infty$ and $m^*((L_0, L_1)) < +\infty$, is found in appendix 1, by S. Kotani, to~\cite{KoWa}. \qed 
 
\subsection{The string associated with a CMRW (definition~\ref{defctgd})} This subsection shows how the holding intensities for the sites are related to the string measure for the continuous time martingale Markov random walk on a finite state space. Let $S = \{i_1,\ldots, i_{M}\}
\subset \mathbb{R}$, $i_1 < \ldots < i_{M}$. Assume there are two constants $0 < c < C < +\infty$ where the inequalities are strict such that $c < \min_{j \in \{2, \ldots, M-1\}} \epsilon_j \leq \max_{j \in \{2, \ldots, M-1\}} \epsilon_j < C$ and that $\min_{j \in \{2, \ldots, M\}} i_j - i_{j-1} > 0$.  Let 

\begin{equation}\label{eqmstep} \left. \begin{array}{l} \tilde{{\cal M}}(x) = \left\{\begin{array}{ll} 0  & x < i_2 \\ 
 \sum_{j=2}^k \epsilon_j  & i_k \leq x < i_{k+1}
\quad k = 2,\ldots, M-2  \\
\sum_{j=2}^{M-1} \epsilon_j &  x \geq i_{M-1} 
\end{array} \right. \\
{\cal M}(x) = \tilde{\cal M}(x) - \tilde{\cal M}(0). \end{array}\right. 
\end{equation}

\noindent In the notation of the definition, $l_0 = i_2$ and $l_1 = i_{M-1}$. Let $L_1 = i_M > i_{M-1}$ and $L_0 = i_1 < i_2$. Let $m^*$ denote the string defined from ${\cal M}$ by equation (\ref{eqmfromM2}). Then  $m^*(\{i_j\}) = \epsilon_j$ for $j = 2, \ldots, M-1$, $m^*(\{i_1\}) = m^*(\{i_M\}) = +\infty$. The generator ${\cal G}^*$ satisfies 

\begin{equation}\label{eqpreop} \left\{\begin{array}{ll}   {\cal G}^*f(i_j) =  
 \frac{1}{m^*(\{i_j\})} \left(\frac{f(i_{j+1}) - f(i_j)}{i_{j+1}
- i_j} - \frac{f(i_j) - f(i_{j-1})}{i_j - i_{j-1}}\right) & j \in \{2,
\ldots, M-1\}  \\ {\cal G}^* f(i_1) = {\cal G}^* f(i_M) = 0 & \end{array}\right.
  \end{equation}

\noindent  This can be expressed as

\begin{equation}\label{eqgexp} {\cal G}^* f(i_j) = \frac{1}{ 
m^*(\{ i_j \})}\frac{(i_{j+1} - i_{j-1})}{(i_j - i_{j-1}) (i_{j+1} - i_j)}
\left(\frac{i_j - i_{j-1}}{i_{j+1} - i_{j-1}}f(i_{j+1}) - f(i_j) + \frac{i_{j+1}
- i_j}{i_{j+1} - i_{j-1}}f(i_{j-1})\right) \end{equation}

\noindent The operator ${\cal G}^*$ is the generator of a CMRW (definition~\ref{defctgd}) with holding intensities 

\begin{equation}\label{eqholding} \left\{ \begin{array}{ll} \lambda(i_j) : =  \frac{(i_{j+1} - i_{j-1})}{ 
m(\{ i_j \} )(i_j - i_{j-1})(i_{j+1} - i_j)} & j = 2,\ldots, M-1 \\
\lambda(i_1) = \lambda(i_M) = 0 & \end{array}\right. 
\end{equation} 

\noindent For $j \in \{ 2, \ldots, M-1\}$, when it jumps from $i_j$, it jumps to
$i_{j-1}$ with probability $\frac{i_{j+1} - i_j}{i_{j+1} - i_{j-1}}$ and to
$i_{j+1}$ with probability $\frac{i_j - i_{j-1}}{i_{j+1} - i_{j-1}}$. Note that
\[ \frac{i_{j+1} - i_j}{i_{j+1} - i_{j-1}}i_{j-1} + \frac{i_j - i_{j-1}}{i_{j+1}
- i_{j-1}}i_{j+1} = i_j \qquad j = 2, \ldots, M-1.\]
\noindent If the process arrives at $i_M$, it remains there for all subsequent time; if it arrives at $i_1$, it remains there for all subsequent time. The process is therefore a martingale.  
 
\subsection{The transition semigroup and differential equations} Let $X$ denote the gap diffusion, definition~\ref{defmpsting}, with state space $(-\infty, +\infty) \cap [L_0, L_1]$, associated with the string measure $m^*$. Let ${\cal P}$ denote the transition probability function; that is, ${\cal P}(s; x, A) = \mathbb{P}(X_s \in A | X_0 = x)$.  Let ${\cal G}^*$  denote the second order differential operator  corresponding to $m^*$.  

Consider $L_1 - L_0 < +\infty$, $m^*((L_0, L_1)) < +\infty$.   For $A \in {\cal B}((L_0,L_1))$, ${\cal P}  (s;.,A) \in L^2([L_0, L_1], m^*)$, since ${\cal P}(s;L_0,A) = {\cal P}(s;L_1,A) = 0$. Let $c_k(s,A) = \langle {\cal P} (s;.,A), \eta_k \rangle_{L^2([l_0, L_1], m^*)}$ and let $c_k(A) = c_k(0,A) = \int_A \eta_k(y) m^*(dy)$. 

Since ${\cal G}^*$ is the infinitesimal generator of $X$, it follows that $\frac{\partial}{\partial s} {\cal P} (s;x,A) = {\cal G}^*{\cal P} (s;x,A)$, from which 
\begin{equation}\label{eqceee} \left\{ \begin{array}{l} \frac{\partial}{\partial s} c_k(s,A) = -\lambda_k c_k(s,A)\\ c_k(0,A) = \int \eta_k(x) m(dx) \end{array}\right. \end{equation}  and hence

\[ {\cal P}  (s;x,A) = \sum_{k = 1}^\infty c_k(A)e^{- \lambda_k t} \eta_k(x).\]

\begin{Lmm} For any $A \subset (L_0, L_1)$, for all $s > 0$, $ {\cal P}(s;.,A) \in {\cal D}({\cal G}^*)$.
\end{Lmm}

\paragraph{Proof} For $L_1 - L_0 < +\infty$, from the definition, the only point that has to be proved is that $\int_{L_0-}^{L_1} ({\cal G}^*  {\cal P}(s;x,A))^2 m^*(dx) < +\infty$. This follows, because ${\cal G}^*{\cal P} (s;x,A) = -\sum_{k=1}^\infty c_k(A)\lambda_k e^{-\lambda_k t} \eta_k(x)$. Since $\lambda^2 e^{-2\lambda s} \leq \frac{1}{s^2}e^{-2}$ (the bound independent of $\lambda$, obtained by differentiation), it follows that 
\[ \int_{L_0-}^{L_1} ({\cal G}^* {\cal P}(s;x,A))^2 m^*(dx) < \frac{1}{s^2}e^{-2}  \int_{L_0-}^{L_1} ( {\cal P}(s;x,A))^2 m^*(dx) < \frac{1}{s^2}e^{-2}m^*((L_0,L_1))\] 

\noindent and the result follows directly from the definitions for $L_1 - L_0 < +\infty$. \vspace{5mm}

\noindent For $L_1 - L_0 = \infty$, the important feature is the upper bound $\lambda^2 e^{-2\lambda s} < \frac{1}{s^2} e^{-2}$. A limiting sequence may be employed; for $-\infty < z_1 < z_2 < +\infty$, $L_0 \leq z_1 < z_2 \leq L_1$ and $A \in {\cal B}((L_0,L_1))$,  let $X$ denote the process generated by ${\cal G}^*$, let $\tau_{z_1} = \inf\{s | X_s \leq z_1\}$, $\tau_{z_2} = \inf\{ s | X_s \geq z_2\}$ and let $\tilde{X}^{z_1,z_2}$ be the process defined by

\begin{equation}\label{eqprocredef} \tilde{X}^{z_1,z_2}_s = \left\{ \begin{array}{ll} X_s & 0 \leq s < \tau_{z_1} \wedge \tau_{z_2} \\ z_1 & s \geq \tau_{z_1} \wedge \tau_{z_2} \qquad \tau_{z_1} \leq \tau_{z_2} \\ z_2 & s \geq \tau_{z_1} \wedge \tau_{z_2} \qquad \tau_{z_1} > \tau_{z_2}. \end{array}\right. \end{equation}

\noindent Then $\tilde{X}^{z_1,z_2}$ is the process associated with the string measure

\begin{equation}\label{eqstringrec} m_{z_1,z_2}^*(A) = \left\{ \begin{array}{ll} m^*((z_1,z_2) \cap A) & z_1 \not \in A, \quad z_2 \not \in A\\ +\infty & z_1 \in A\quad \mbox{or} \quad z_2 \in A. \end{array}\right. \end{equation}

\noindent Set 

\[ \tilde{{\cal P}}_{z_1,z_2}(s;x,A) = \mathbb{P}\left ( \{X_s \in A \cap (z_1,z_2)\} \cap \left \{z_1 < \inf_{0 \leq r \leq s} X_r \leq \sup_{0 \leq r \leq s} X_r < z_2\right \}| \tilde{X}_0 = x \right ).\]

\noindent Let ${\cal P}_{z_1,z_2}$ denote the transition probability for the process associated with the string measure $m_{z_1,z_2}^*$ and note that $\tilde{{\cal P}}_{z_1,z_2}(s;x,A) = {\cal P}_{z_1,z_2}(s;x,A\cap(z_1,z_2))$. The quantity $\tilde{{\cal P}}_{z_1,z_2}(s,x,A)$ is non decreasing as $z_1 \downarrow L_0$ and $z_1 \uparrow L_1$ and \[{\cal P}(s;x,A) = \lim_{z_1 \downarrow L_0} \lim_{z_2 \uparrow L_1} \tilde{{\cal P}}_{z_1,z_2}(s;x;A) \quad \forall x \in (L_0,L_1), \quad \forall A \in {\cal B}((L_0,L_1)).\] The result follows from the definition of ${\cal D}({\cal G}^*)$ together with the upper bound on $\lambda^2 e^{-2\lambda s}$ which is independent of the eigenvalue for $s > 0$. \qed \vspace{5mm}

\noindent For $x \leq L_0$, ${\cal P}(s, x, A) = 1$ if $L_0 \in A$ and $0$ otherwise. For $x \geq L_1$, ${\cal P}(s, x, A) = 1$ if $L_1 \in A$ and $0$ otherwise; the process $X$ is a martingale if the initial condition $x$ satisfies $x \in [L_0, L_1]$. \vspace{5mm}

\noindent The following lemma helps to describe the structure of the transition probability; equation (\ref{eqpdef}) shows that the transition probability may be decomposed as a symmetric function, within the domain of the operator in both space variables for all $t > 0$, integrated against the measure. This decomposition is well established (see, for example, Borodin and Salminen~\cite{BS} for the absolutely continuous setting).

\begin{Lmm} \label{lmmq} There exists a function $q$ such that $q(s;x,y) = q(s;y,x)$ for all $(x,y)$, $s \geq 0$, satisfying $q(s;.,.) \in {\cal D}({\cal G}^*) \times {\cal D}({\cal G}^*)$ for each $s > 0$, such that, for $A \in {\cal B}((L_0,L_1))$, the transition probability may be written 

\begin{equation}\label{eqpdef} {\cal P}(t;x,A) = \int_A q(t;x,y) m^*(dy) \qquad A \in {\cal B}((L_0,L_1)) \qquad x \in \mathbb{R}.\end{equation}

\noindent For $L_1 - L_0 < +\infty$ and $m^*((L_0,L_1)) < +\infty$, let $(-\lambda_k, \eta_k)_{k \geq 0}$ denote the eigenvalues and eigenfunctions corresponding to the complete orthonormal base of the operator ${\cal G}^*$ associated with $m^*$. The function $q(t;x,y)$ may be written

\begin{equation}\label{eqqdef} q(t;x,y) = \sum_k e^{-\lambda_k t} \eta_k(x) \eta_k(y).\end{equation}
\end{Lmm}
 
\paragraph{Proof} First consider $L_1 - L_0 < +\infty$. From equation (\ref{eqceee}), it follows that   ${\cal P}(t;x,A) = \sum_k c_k(t,A) \eta_k(x)$ where \[c_k(t,A) = e^{-\lambda_k t} \int_A \eta_k(x) m^*(dx).\]

\noindent It follows that ${\cal P}$ satisfies equation (\ref{eqpdef}) with $q$ defined by equation (\ref{eqqdef}). For $L_1 - L_0 = +\infty$, consider the restricted processes of equation (\ref{eqprocredef}) associated with the restricted strings of equation (\ref{eqstringrec}), with ${\cal P}_{z_1,z_2}$ the transition probability for the restricted process. Since ${\cal P}_{z_1,z_2}(s;x,A \cap (z_1,z_2))$ is non decreasing as $z_1 \downarrow -\infty$ and $z_2 \uparrow +\infty$ for all $A \in {\cal B}((L_0,L_1))$, it follows that the functions $q_{z_1,z_2}(t;x,y)$ from equation (\ref{eqqdef}) for the restricted string are non decreasing as $z_1 \downarrow -\infty$ and $z_1 \uparrow +\infty$ and the function $q(t;.,.) = \lim_{z_1 \downarrow -\infty}\lim_{z_2 \uparrow +\infty} q_{z_1,z_2}(t;.,.)$ satisfies the criteria. \qed \vspace{5mm}

\section{Proof of theorem~\ref{thgdcont}}\label{appthgdcont}

This section is devoted to the proof of theorem~\ref{thgdcont}, that for any probability measure $\mu$ over $\mathbb{R}$ such that $\int_{-\infty}^\infty |x| \mu(dx) < +\infty$, and any fixed $t > 0$, there exists a gap diffusion (definition~\ref{defmpsting}) $X$ such that for all $x \in \mathbb{R}$, $\mathbb{P}(X_t \leq x | X_0 = e_0(\mu)) = \mu((-\infty, x])$. 

The result of theorem~\ref{thgdl2} is used as the starting point for the proof.  The measure $\mu$ is approximated by a sequence $\mu^{(n)}$ where for each $n$, as a function of $x$, $\mu^{(n)}((-\infty, x])$ is a step function. For each $\mu^{(n)}$, theorem~\ref{thgdl2} gives a process, for which there is an  associated string measure $m^{*(n)}$. The characterisation of the Markov process associated with a string $m^*$ is given in definition~\ref{defmpsting}. This turns out to be a useful characterisation for showing that convergence of strings implies convergence of  associated processes. That is the subject of lemma~\ref{lmmloctimeconv}. After this is established, the remainder of the proof of theorem~\ref{thgdcont} involves showing that convergence of $\mu^{(n)}$ to a measure $\mu$ implies existence of a convergent subsequence of strings.

\begin{Lmm} \label{lmmloctimeconv} Let $(m^{*(n)})_{n \geq 1}$ denote a sequence of strings, defined according to definition~\ref{defks} equation (\ref{eqmfromM2}), with $L_{0,n}$ and $L_{1,n}$ as left and right cemetaries respectively for the string $m^{*(n)}$,  $L_{0,n} \downarrow L_0 > -\infty$ and  $L_{1,n} \uparrow L_1 < +\infty$. Assume that there is a string measure $m^*$, with left and right cemetaries $L_0$ and $L_1$ such that   

\begin{equation}\label{eqconvhyp} \lim_{n \rightarrow +\infty} \sup_{z} |m^{*(n)}((L_{0,n} , z \wedge L_{1,n})) - m^*((L_0  , z \wedge L_1))| = 0.
\end{equation}

\noindent  Let $X^{(n)}$ denote the Markov process associated with  $m^{(n)*}$ according to definition~\ref{defmpsting} and $X$ the Markov process associated with $m^*$. Then for each $z \in (L_0, L_1)$, any sequence $(z_n)_{n \geq 1}$ such that $z_n \rightarrow z$ and $z_n \in (L_{0,n},L_{1,n})$ for each $n$ and any $t < +\infty$, 

\[ \lim_{n \rightarrow +\infty} \sup_{w \in \mathbb{R}} \left | \mathbb{P}\left(   X^{(n)}_t  \leq w | X^{(n)}_0 = z_n \right ) - \mathbb{P}\left( X_t \leq w |X_0 = z\right ) \right |   = 0.\]

\end{Lmm}

\paragraph{Proof}  Let $W$ denote a standard Wiener process started from $0$ and let $\phi(s,a)$ denote its local time at site $a \in \mathbb{R}$, at time $s \geq 0$. Let 
\[ T^{(n)}(z,s) = \int_{\mathbb{R}}\phi (s,a-z)m^{*(n)}(da), \quad T(z,s) = \int_{\mathbb{R}}\phi (s,a - z)m^{*}(da).\]

\noindent  Note that, for each $a$, $\phi (r,a)$ is nondecreasing in $r$, continuous in both $r$ and $a$ and $\lim_{r\rightarrow +\infty} \phi (r,a) = +\infty$ for each $a \in \mathbb{R}$. Define 
\[ T^{(n)-1}(z,s) = \inf \left \{r |  \int_{\mathbb{R}} \phi (r,a-z) m^{*(n)}(da) \geq s \right \}, \quad  T^{-1}(z,s) = \inf \left \{r |  \int_{\mathbb{R}} \phi (r,a- z) m^{*}(da) \geq s \right \}.\] 

\noindent Then, from definition~\ref{defmpsting}, $X^{(n)}$ has representation $X^{(n)}(t,z) = z + W(T^{(n)-1}(z,t))$, while $X(t,z)$ has representation $X(t,z) = z + W(T^{-1}(t,z))$.  The result now follows from this characterisation of $X^{(n)}$ and $X$, together with continuity of $W$ and the continuity of its local time giving $\lim_{n \rightarrow +\infty} T^{(n)-1}(z_n,t) \rightarrow T^{-1}(z,t)$ almost surely, from which convergence of the laws follows. \qed

\paragraph{Proof of theorem~\ref{thgdcont}} Let $\mu$ denote the target measure. Define the points $x_{n,j}$ as

\[ \left\{\begin{array}{l}  x_{n,1} = \inf\{z | \mu((-\infty,z]) \geq \frac{1}{2^n}\}    \\
x_{n,j+1} =  \inf \left \{ z > x_{n,j} | \mu((-\infty, z]) -
\mu((-\infty, x_{n,j}]) \geq \frac{1}{2^n} \right \}\\ \hspace{10mm} \mbox{until either} \; \mu((-\infty, x_{n,j+1}]) = 1 \; \mbox{or} \;  1 - \mu ((-\infty, x_{n,j+1}]) < \frac{1}{2^n}\\
 M_n = \max\{k : \mu((-\infty, x_{n,k}]) < 1\}.  \end{array} \right. \]

\noindent Let $\underline{p}^{(n)} = (p^{(n)}_1, \ldots, p^{(n)}_{M_n})$ denote the probability mass function defined by 

\[\left\{ \begin{array}{l} p^{(n)}_1= \mu((-\infty, x_{n,1}]) \\ p^{(n)}_j =   \mu((x_{n,j-1}, x_{n,j}]) \qquad j = 2, \ldots, M_n - 1  
  \\ p^{(n)}_{M_n} = 1 - \mu((-\infty, x_{n,M_{n-1}}]), \end{array}\right. \]
  \noindent over $\{x_{n,1}, \ldots, x_{n,M_n}\}$ with expectation 
\[e_0(\underline{p}^{(n)}) = \sum_{j = 1}^{M_n} x_{n,j}p^{(n)}_j.\]
                      
\noindent Then, by theorem~\ref{thgdl2},  for each $n$ there exists a CMRW (definition~\ref{defctgd}) $X^{(n)}$ with a finite state space $S_n = \{x_{n,1} \ldots, x_{n,M_n}
\}$ 
such that $\mathbb{P}(X^{(n)}_t = x_{n,j} | X^{(n)}_0 = e_0(\underline{p}^{(n)})) =
p^{(n)}_j$. 

\noindent Note that, for $j = 1, \ldots, M_{n-1}$, 

\[ \mathbb{P} \left (X_t^{(n)} \leq x_{n,j} | X_0^{(n)} = e_0(\underline{p}^{(n)}) \right ) = \mu ((-\infty,x_{n,j}]) \qquad j = 1,\ldots, M_n.\]

\noindent Furthermore, it follows that 

\[e_0(\underline{p}^{(n)}) \rightarrow e_0(\mu) := \int_{-\infty}^{\infty} x \mu(dx).\]

\noindent Note that, by construction, for any $n < + \infty$, $- \infty < x_{n,1} \leq x_{n,M_n} < +\infty$, $x_{n,1}$ is non-inreasing in $n$ and $x_{n,M_n}$ is non-decreasing in $n$. Let   $Q_0 = \lim_{n \rightarrow +\infty} x_{n,1}$ and  $Q_1 = \lim_{n \rightarrow +\infty} x_{n,M_n}$; $Q_0$ possibly $-\infty$ and  $Q_1$ possibly $+\infty$. \vspace{5mm}

\noindent Let 

\begin{equation}\label{eqjumpdef} {\cal S}  = \{x | \mu(\{x\}) > 0\}.\end{equation} 

\noindent  ${\cal S}$ is the set of atoms of $\mu$. Let $m^{*(n)}$ denote the string corresponding to the process $X^{(n)}$. Then there is a non decreasing function ${\cal M}^{(n)}$, with ${\cal M}^{(n)}(0 ) = 0$, ${\cal M}^{(n)}(x_{n,2}) > -\infty$, ${\cal M}^{(n)}(x_{n, M_{n-1}}) < +\infty$ whose points of increase are  $\{x_{n,2}, \ldots, x_{n,M_n-1}\}$, from which the string $m^{*(n)}$ may be constructed using equation (\ref{eqmfromM2}), with $L_0 = x_{n,1}$ and  $L_1 = x_{n,M_n}$.

\begin{Lmm}\label{lmmceeup} For any $-\infty < x < y < +\infty$,  there is a constant $C(x,y) < +\infty$ such that 

\begin{equation}\label{eqceeup} \sup_{n \geq 1} m^{*(n)}((x \vee x_{n,1}, y \wedge x_{n,M_n})) < C(x,y) <  +\infty
 \end{equation}
 \end{Lmm}
 
\paragraph{Proof of lemma~\ref{lmmceeup}} Assume not, then   there exists a pair $(x,z)$ such that  $-\infty < x < z < +\infty$ and either an $n \geq 1$ such that $m^{*(n)}((x \vee x_{n,2} ,z \wedge x_{n,M_n})) = +\infty$ or a subsequence $(n_j)_{j \geq 1}$ such that $\lim_{j \rightarrow +\infty} m^{*(n_j)}((x \vee x_{n_j,1}, z \wedge x_{n_j,M_{n_j}})) = +\infty$.\vspace{5mm}

\noindent Assume the first case. Then there is an $n > 1$ and a point $x_{n,k} \in (x,z)$, where $2 \leq k \leq M_n - 1$, such that $m^{*(n)}(\{ x_{n,k} \} )  = +\infty$. If $x_{n,k} > e_0(\underline{p}^{(n)})$, then from equation (\ref{eqholding}), it follows that \[ \mathbb{P} \left (X^{(n)}_t \leq z | X_0^{(n)} = e_0(\underline{p}^{(n)})\right ) = 1 \qquad \forall t \geq 0,\]

\noindent and, consequently $\mu((-\infty, z]) = 1$ for some $z \leq x_{n,M_{n} - 1}$,
contradicting the hypotheses on $\mu$. 

If $x_{n,k} < e_0(\underline{p}^{(n)})$, then $\mathbb{P} \left (X_t \geq x_{n,k}|X_0 = e_0(\underline{p}^{(n)}) \right ) = 1$ for $x_{n,k} \geq x_{n,2}$, again a contradiction. \vspace{5mm}

\noindent Now consider the second case; suppose that for all $(x^\prime , z^\prime )$ such that $-\infty < x^\prime < z^\prime  < +\infty$, $m^{*(n)}((x^\prime \vee x_{n,1}, z^\prime \wedge x_{n,M_n})) < +\infty$ for each $n$, but that there is a pair $(x^\prime,z^\prime)$, $x^\prime > Q_0$ and $z^\prime < Q_1$ and a subsequence $(n_j)_{j \geq 1}$  such that $\lim_{j \rightarrow +\infty} m^{*(n_j)}((x^\prime \vee x_{n_j,1}, z^\prime \wedge x_{n_j,M_{n_j}})) = \infty$. Let 

\[\tilde{z} = \inf\{z^\prime | \lim_{j \rightarrow +\infty} m^{*(n_j)}((x^\prime \vee x_{n_j,1}, z^\prime \wedge x_{n_j,M_{n_j}})) = \infty\}\]

\noindent and, for this $\tilde{z}$, let \[\tilde{x} = \lim_{\epsilon \downarrow 0} \left(\sup\{x^\prime  | \lim_{j \rightarrow +\infty} m^{*(n_j)}((x^\prime \vee x_{n_j,1}, (\tilde{z} + \epsilon) \wedge x_{n_j,M_{n_j}})) = \infty\}\right ).\]

\noindent Note that $\tilde{x} = \tilde{z}$. Consider the representation $X^{(n)}(s,x) = W(T^{(n)-1}(x,s),x)$ from definition~\ref{defmpsting}, where $T^{-1}$ is the inverse function of $T^{(n)}(x,s) = \int_{\mathbb{R}} \phi^{(x)}(s,z)m^{*(n)}(dz)$. By the continuity of $\phi^{(x)}$, it follows that $\lim_{j \rightarrow +\infty} T^{(n_j)}(x,s) = +\infty$ if $\phi^{(x)}(s,\tilde{z}) > 0$. If $\tilde{z} > e_0(\mu)$, it follows that $\lim_{j \rightarrow +\infty}\mathbb{P}(X^{(n_j)}_t > \tilde{z} | X_0^{(n_j)} = e_0(\underline{p}^{(n_j)} )) = 0$ and if $\tilde{z} < e_0(\mu)$ then  $\lim_{j \rightarrow +\infty}\mathbb{P}(X^{(n_j)}_t < \tilde{z} | X_0^{(n_j)} = e_0(\underline{p}^{(n_j)} )) = 0$, a contradiction in both cases.
 
  Lemma~\ref{lmmceeup} is proved. \qed  

\paragraph{Proof of theorem~\ref{thgdcont} (continued)}  Recall the definition of ${\cal S}$ in equation (\ref{eqjumpdef}), the set of points where $\mu$ has an atom.  Since $\mu$ is a probability measure,
there are at most a countable number of elements of ${\cal S}$. Note also that
for each $x \in {\cal S}$, there is an $N$ such that for all $n > N$, $x =
x_{n,k_n(x)}$ for some $k_n(x)$.

The remainder of the proof splits the string functions $({\cal M}^{(n)})_{n \geq 1}$ into two parts, a part that converges by the Ascoli Arzela theorem to a continuous function and a part that converges to a function that increases by  jumps.  

Set $x_{n,0} = x_{n,1} - \frac{1}{2^n}$ and define $\tilde{{\cal M}}^{(n)}$ by $\tilde{{\cal M}}^{(n)}(x_{n,0}) = 0$ and, for $j =  0, 1, \ldots, M_n - 1$,  

\[ \tilde{{\cal M}}^{(n)}(x) = \left\{\begin{array}{ll} \tilde{{\cal M}}^{(n)}(x_{n,j}) +
\frac{x - x_{n,j}}{x_{n,j+1} - x_{n,j}}({\cal M}^{(n)}(x_{n,j+1}) - {\cal M}^{(n)}(x_{n,j})) &
x \in [x_{n,j}, x_{n,j+1}], \quad  x_{n,j+1} \not \in {\cal S} \\
\tilde{{\cal M}}^{(n)}(x_{n,j}) & x \in [x_{n,j}, x_{n,j+1}], \quad x_{n,j+1} \in S.
\end{array}\right. \]

\noindent  Let $\widehat{\tilde{{\cal M}}}^{(n)}(x) = \tilde{{\cal M}}^{(n)}(x) -  \tilde{{\cal M}}^{(n)}(0)$.  For $-\infty < x < y - \delta < +\infty$, let  

\[ \epsilon_{x,y} (\delta) = \limsup_{n \rightarrow +\infty} \sup_{z \in [x_{n,1} \vee x , x_{n,M_n} \wedge y - \delta)}\sup_{z  \leq z^\prime <
z+\delta}|\tilde{{\cal M}}^{(n)}(z^\prime) - \tilde{{\cal M}}^{(n)}(z)|.\]

\begin{Lmm}\label{lmmepszero} \[ \lim_{\delta \rightarrow 0} \epsilon_{x,y}(\delta) = 0 \qquad \forall -\infty < x < y < +\infty.\]
\end{Lmm}

\paragraph{Proof of lemma~\ref{lmmepszero}} If $\limsup_{\delta \rightarrow 0} \epsilon_{x,y} (\delta) = c > 0$, for some $c > 0$ and some $(x,y)$, then there is a
sequence $(x_{n_j}, \delta_{n_j}, \tilde{{\cal M}}^{(n_j)})_{j \geq 1}$ such that $\delta_{n_j}
\downarrow 0$ and such that there exists an $x$ and $y$ such that 

\[ \lim_{j \rightarrow +\infty} (\tilde{{\cal M}}^{(n_j)}(x_{n_j} + \delta_{n_j}) -
\tilde{{\cal M}}^{(n_j)}(x_{n_j})) \geq c > 0, \qquad -\infty < x < \inf_j x_{n_j} \leq  \sup_j x_{n_j} < y < +\infty. \]

\noindent  Consider a limit point $x^*$ of
$x_{n_j}$. If $x^* \not \in S$, then 

\[ \lim_{n \rightarrow +\infty} \mathbb{P} \left(X_t^{(n)} = x^* | X_0 = e_0\left(\underline{p}^{(n)}\right)\right) = 0.\] 

\noindent If $x^* \in {\cal S}$ then, by the construction of $\underline{p}^{(n)}$, there exists an $N$ such that for each $N \geq n$,

\begin{equation}\label{eqnolim} \left | \mathbb{P}\left (X^{(n)}_t = x^*|X_0 = e_0 \left (\underline{p}^{(n)} \right ) \right ) -  \mu (\{x^*\} ) \right | \leq \frac{1}{2^{n-1}}.
\end{equation}

\noindent Note that 
\begin{eqnarray*}\lefteqn{ \mathbb{P} \left(X_t^{(n)} \in [x^* - \delta, x^* + \delta] |X_0 = e_0 \left (\underline{p}^{(n)} \right )\right)}\\&& = \mathbb{P} \left(X_t^{(n)} = x^* |X_0^{(n)} = e_0 \left (\underline{p}^{(n)} \right )\right) + \mathbb{P} \left(X_t^{(n)} \in [x^* - \delta, x^* + \delta]\backslash\{x^*\} |X_0^{(n)} = e_0 \left (\underline{p}^{(n)} \right )\right) \end{eqnarray*}

\noindent and that 

\[ \mathbb{P} \left(X_t^{(n)} = x^* |X_0^{(n)} = e_0 \left (\underline{p}^{(n)} \right )\right) \rightarrow \mu (\{x^*\}). \]

\noindent It follows that, in all cases,   

\begin{equation}\label{eqstineq} \lim_{\delta \downarrow 0} \lim_{n \rightarrow +\infty}\mathbb{P}\left (X^{(n)}_t \in
[x^* -\delta, x^* +\delta] \backslash\{x^*\} |X_0^{(n)} = e_0 \left (\underline{p}^{(n)} \right )\right) = 0.
\end{equation}
 
\noindent The aim is to show that if $\lim_{\delta \rightarrow 0}\epsilon_{x,y}(\delta) = c > 0$, then equation (\ref{eqstineq}) does not hold, leading to a contradiction.

\noindent Set 

\[ g^{(n)} (s,x) = \mathbb{P} \left (X_s^{(n)} \in [x^* - \delta, x^* + \delta] \backslash \{x^*\}   |X_0^{(n)} = x \right).\]

\noindent   Let  $q^{(n)}(t;x,y)$ be the element of ${\cal D} ({\cal G}^{*(n)}) \times {\cal D} ({\cal G}^{*(n)})$ from lemma~\ref{lmmq} such that ${\cal P}^{(n)} (t;x,A) = \int_A q^{(n)} (t;x,y)m^{*(n)}(dy)$ for $A \in {\cal B}((x_{n,1},x_{n,M_n}))$. It follows that 

\[ g^{(n)}(s,x) = \sum_{x_{n,k} \in [x^* - \delta, x^* +\delta] \backslash \{x^*\}} q^{(n)}(x, x_{n,k}) m^{*(n)}(\{x_{n,k}\}).\]    

\noindent If $\lim_{\delta \rightarrow 0} \epsilon_{z_1,z_2}(\delta) = c > 0$, then for any $x \in (z_1,z_2)$, 

\begin{equation}\label{eqstrlb} \lim_{\delta \downarrow 0} \liminf_{n \rightarrow +\infty} g^{(n)}(t,x) \geq c   \liminf_{n \rightarrow +\infty} q^{(n)}(t;x,x^*).\end{equation}

\noindent It holds that 
\[ \liminf_{n \rightarrow +\infty} q^{(n)}(t; e_0(\underline{p}^{(n)}),x^*) = \liminf_{n \rightarrow +\infty} q^{(n)}(t; e_0(\mu), x^*).\]

\noindent The following argument shows that $\liminf_{n \rightarrow +\infty} q^{(n)}(t; e_0(\mu), x^*) > 0$, where the inequality is strict. If $e_0(\mu) < x^*$, let $\tau^{(n)} = \inf\{s | X^{(n)}_s \geq x^*\}$ and $y_n^* = \inf\{x_{n,k}|x_{n,k} \geq x^* \}$. If $e_0(\mu) > x^*$, let $\tau^{(n)} = \inf\{s | X^{(n)}_s \leq x^*\}$ and $y^*_n = \inf\{x_{n,k} | x_{n,k} \leq x^*\}$. Then

\[ q^{(n)}(t;e_0(\mu),y^*_n) = \int_0^t \mathbb{P}(\tau^{(n)} \in ds|X^{(n)}_0 = e_0(\mu))q^{(n)}(t-s;y^*_n,y^*_n).\]

\noindent It follows the construction that $\liminf_{n \rightarrow +\infty} q^{(n)}(s;y^*_n, y^*_n) > 0$; otherwise it follows from the construction that $\liminf_{n \rightarrow +\infty} q^{(n)}(s;y^*_n,.) \equiv 0$ and, using the fact that $m^{*(n)}((x_{n,1}\vee a, x_{n,M_n} \wedge b)) < C(a,b) < +\infty$, that $\liminf_{n \rightarrow +\infty} {\cal P}^{(n)}(s;y^*_n, (a \vee x_{n,1}, b \wedge x_{n,M_n})) = 0$ for all $-\infty < a < b < +\infty$, $a \geq Q_0$ and $b \leq Q_1$, implying that there is a subsequence $(n_j)_{j \geq 1}$ such that $\lim_{j \rightarrow +\infty} {\cal P}^{(n_j)}(s;y^*_{n_j}, (Q_0,Q_1)) = 0$, from which it follows that $\mu(\{Q_0\}) + \mu(\{Q_1\}) = 1$, from which $m^{*(n)}((Q_0,Q_1)) = 0$ for all $n \geq 1$, which is a contradiction. 

It follows that $\lim_{\delta \rightarrow 0}   \epsilon_{x,y}(\delta) = 0$ for all $-\infty < x < y < +\infty$ and the lemma is proved.   \qed

\paragraph{Proof of theorem~\ref{thgdcont} (concluded)} Recall that $Q_0 = \lim_{n \rightarrow +\infty} x_{n,1}$ and $Q_1 = \lim_{n \rightarrow +\infty} x_{n,M_n}$. It follows, by the Arzelo-Ascoli theorem that for any $-\infty < x < y < +\infty$ satisfying $x \geq Q_0$ and $y \leq Q_1$, there is a subsequence
$\widehat{\tilde{\cal{M}}}^{(n_k)}$ and a limit point $\widehat{\tilde{\cal{M}}}_{x,y}$ such that 

\[ \lim_{k \rightarrow +\infty} \sup_{z \in [x,y]} \left |\widehat{\tilde{\cal{M}}}_{x,y}(z) - \widehat{\tilde{\cal{M}}}^{(n_k)}(z) \right | = 0.\]

\noindent  For $-\infty < x < y < +\infty$, it follows from weak convergence and compactness arguments that there exists a subsequence $\{m^{*(n_{k_j})}(\{z\})|z \in {\cal S}\}$, where $(n_{k_j})_{j \geq 1}$ is a subsequence of the $(n_j)_{j \geq 1}$ used for the continuous part, 
and a limit point $\{\theta_{x,y} (z) | z \in {\cal S}\}$ such that  for any $\beta$ satisfying $|\beta(z)|  < \frac{1}{1 + C(0,0 \vee z)^2 + C(0 \wedge z,0)^2}$, 
\[ \sum_{z \in {\cal S}: x_{n_{k_j},2} \leq z \leq x_{n_{k_j},M_{n_{k_j}}-1}}
m^{*(n_{k_j})}(\{z\})\beta(z) \stackrel{j \rightarrow +\infty}{\longrightarrow}
\sum_{z \in {\cal S}| Q_0 < z < Q_1} \theta_{x,y} ( z) \beta(z).\]

\noindent Now let 

\[ f_{x,y}(z) = \left\{\begin{array}{ll} \sum_{w \in [0,z] \cap {\cal S}} \theta_{x,y}(w) & z \geq 0 \\ - \sum_{w \in [z , 0) \cap {\cal S}} \theta_{x,y} (w) & z \leq 0 \end{array}\right. \] 

\noindent Define ${\cal M}_{x,y}$ as 

\[  {\cal M}_{x,y}(z) = \left\{ \begin{array}{ll} \widehat{\tilde{\cal M}}_{x,y}(x)  +   f_{x,y}(x) & z \leq x \\ \widehat{\tilde{\cal M}}_{x,y}(z)  +   f_{x,y}(z) & z \in (x,y) \\ \widehat{\tilde{\cal M}}_{x,y}(y)  +   f_{x,y}(y) & z \geq y. \end{array}\right. \]

\noindent  Let ${\cal M}$ be a limit point of ${\cal M}_{x,y}$ in the sense that there exists a sequence $(x_i,y_i)_{i \geq 1}$ such that $x_i \downarrow -\infty$, $y_i \uparrow +\infty$ and $\lim_{i \rightarrow +\infty} \int\gamma(z) ({\cal M}_{x_i,y_i}(z) - {\cal M}(z)) dz = 0$ for all $\gamma$ such that $\int |\gamma(z)|(1 + C(z\wedge 0,0)^2 + C(0, z\vee0)^2 )dz < +\infty$.  Such a limit exists by the usual convexity and compactness arguments. 

Let $m^*$ denote the string measure constructed from ${\cal M}$ using equation (\ref{eqmfromM2}) with $L_0 = Q_0$ and $L_1 = Q_1$. Then for each $-\infty < x < y < +\infty$ such that $Q_0 \leq x < y \leq Q_1$, by construction, there is a subsequence $(m^{*(n_{j})})_{j \geq 1}$ such that   

\begin{equation}\label{eqconvnice} \lim_{j \rightarrow +\infty} \sup_{z \in [x,y]} \left |m^{*(n_{j})}((x \vee x_{n_{j},1},z \wedge x_{n_{j}, M_{n_{j}}}) ) - m^*((x \vee Q_0,z \wedge Q_1)) \right |=
0.\end{equation}

\noindent Let $m_{x,y}^{*(n)}$ be defined as the restriction of $m^{*(n)}$ according to equation (\ref{eqstringrec}) and $X^{(n);x,y}$ the corresponding process according to equation (\ref{eqprocredef}); similarly for $m^*_{x,y}$ and $X^{x,y}$. 
  Let $(n_{j})_{j \geq 1}$ denote a sequence for which equation (\ref{eqconvnice}) holds. It follows from lemma~\ref{lmmloctimeconv} that for any $z \in (x \vee Q_0, y \wedge Q_1)$ and sequence $z_{n_j}$ such that $z_{n_j} \in (x \vee Q_0, y \wedge Q_1)$ for each $j$ and $z_{n,j} \rightarrow z$,  

\[ \lim_{j \rightarrow +\infty} \sup_w \left | \mathbb{P} \left ( X^{(n_j);x,y}_t \leq w|X^{(n_j);x,y}_0 = z_{n_j} \right ) - \mathbb{P}\left (X^{x,y}_t  \leq w |X^{x,y}_0 = z \right ) \right | = 0.\]

\noindent It also follows from the above description that, for any $z_1, z_2 \in (Q_0, Q_1)$,

\[ \lim_{x \downarrow -\infty} \lim_{y \uparrow +\infty} \mathbb{P}(X^{x,y}_t  \leq z_2 | X^{x,y}_0 = z_1) = \mathbb{P} (X_t \leq z_2 | X_0 = z_1).\]

\noindent and for any $z_1,z_2 \in (x_{n,1}, x_{n,M_n})$ that

\[\lim_{x \downarrow -\infty} \lim_{y \uparrow +\infty} \mathbb{P}(X^{(n);x,y}_t \leq z_2|X^{(n);x,y}_0 = z_1) = \mathbb{P} (X^{(n)}_t \leq z_2 | X^{(n)}_0 = z_1).\]

\noindent From this, it follows that for any $z \in \mathbb{R}$ and sequence  $(x_j,y_j)$ such that $x_j \downarrow -\infty$ and $y_j \uparrow +\infty$, there exists a subseqence $(n_j)_{j \geq 0}$ such that 

\begin{eqnarray*}\mathbb{P}(X_t \leq z | X_0 = e_0(\mu)) &=& \lim_{j \rightarrow +\infty} \mathbb{P}(X^{(n_j);x_j,y_j}_t \leq z | X^{(n_j);x_j,y_j}_0 = e_0(\underline{p}^{(n_j)}))\\ &=&  \lim_{j \rightarrow +\infty} \mathbb{P}(X^{(n_j)}_t \leq z | X_0 = e_0(\underline{p}^{(n_j)})) = \mu((-\infty ,z]).
\end{eqnarray*}

\noindent and hence there exists a string measure $m^*$ and a gap diffusion $X$ with infinitesimal generator $\frac{d^2}{dm^* dx}$ that satisfies

\[ \mathbb{P}(X_t \leq z | X_0 = e_0(\mu)) = \mu((-\infty, z]) \qquad \forall z \in \mathbb{R}.\] \qed

 \section{Conclusion} For a given marginal distribution $\mu$ at a given fixed
time $t > 0$, with a well defined expectation $e_0(\mu)$, this article proves existence of a gap diffusion which, when started at $e_0$ at time $0$, has this
distribution at time $t$. The motivation for the problem that is currently receiving substantial attention, is from mathematical finance. It is described in the introduction; a list of call option prices for a given maturity $t$ implies, at least approximately, a probability distribution for the stock price at maturity. For purposes of calibration, it is of interest to
recover a time homogeneous equivalent martingale measure that produces this
marginal distribution at the fixed time in question.
 
The article gives existence in full generality for measures with compact support. The theorem of interest to practitioners is theorem~\ref{thgdl2}, on a finite state space.  Although the theorem only states existence, the crux of the problem is locating a vector $\underline{\lambda}$ that satisfies equation (\ref{eqfpeqn}). For fixed $r$, this is a system of polynomial equations, each of degree $r$ and algorithms exist for locating solutions. For fixed $r$, a solution provides a martingale that has the correct marginal for a time $T_r$ with $\Gamma(r, \frac{t}{r})$ distribution. As $r \rightarrow +\infty$, the random variable $T_r$ approaches the deterministic time $t$; for all $\epsilon > 0$, $\lim_{r \rightarrow +\infty} \mathbb{P}(|T_r - t|> \epsilon) = 0$. 

It is the limit of $\underline{\rho}  = r\underline{\lambda}^{(r)}$ that is of interest and hence a limit point, as $r \rightarrow +\infty$ of solutions to

\[ \underline{\rho}  = r{\cal F} \left(\underline{p}N^{r-1} \left ( \frac{\underline{\rho}}{r} \right ) \right ).\]

\noindent  Suitable approximations to this problem should be achievable, since the structure of $N$ is relatively straightforward. It is therefore possible  that the ideas used here in the proof of existence may provide a method for estimating an appropriate martingale diffusion.

\paragraph{Acknowledgements} The author acknowledges the seminar that Jan Ob\l{}\'{o}j  gave in Warsaw on 24th March 2011, which introduced the author to the problem and thanks him for an interesting presentation. He also thanks him for a communication by  electronic mail on 19th May 2011, informing him about the article by Ekström, Hobson, Janson and Tysk~\cite{EHJT} and the article by Monroe~\cite{Mon}.

The author also thanks two anonymous referees, both of whom made invaluable comments, which led to substantial improvements with the revision and, in particular to a referee who pointed out that the method of proof extended easily from measures with compact support to measures over the whole real line and who also suggested using the characterisation of the Markov process in terms of the string and local time of the Wiener process, which proofs that were substantially more elegant than those of the original.

\setlength{\baselineskip}{2ex}

\end{document}